\documentclass{article}

\usepackage{amsmath,amsfonts,amssymb,a4,enumerate,epsfig,theorem,psfrag}
\usepackage[latin1]{inputenc}

\newcommand{\sign}{{\operatorname{sign}}}

\newcommand{\transpose}{{\ast}}

\newcommand{\NN}{{\mathbb{N}}}

\newcommand{\RR}{{\mathbb{R}}}

\newcommand{\Bb}{{\cal{B}}}

\newcommand{\Dd}{{\cal{D}}}
\newcommand{\Ee}{{\cal{E}}}

\newcommand{\Kk}{{\cal{K}}}

\newcommand{\Rr}{{\cal{R}}}

\newcommand{\Tt}{{\cal{T}}}

\newcommand{\Uu}{{\cal{U}}}

\newcommand{\Xx}{{\cal{X}}}
\newcommand{\Yy}{{\cal{Y}}}

\newcommand{\ILa}{\Tt_\Lambda}
\newcommand{\bF}{{\mathbf{F}}}

\newcommand{\bQ}{{\mathbf{Q}}}
\newcommand{\bR}{{\mathbf{R}}}
\newcommand{\bW}{{\mathbf{W}}}
\newcommand{\bc}{{\mathbf{c}}}

\newcommand{\ver}{{a/10}}
\newcommand{\bomega}{\boldsymbol{\omega}}
\newcommand{\bXx}{\boldsymbol{\Xx}}
\newcommand{\bUu}{\boldsymbol{\Uu}}
\newcommand{\bRr}{\boldsymbol{\Rr}}

\newcommand{\xxxx}{\Xx}

\newcommand{\proof}{{\smallskip\noindent\bf Proof: }}

\def\qed{\unskip\nobreak\hfil\penalty50\hskip1.75em\null\nobreak\hfil
$\blacksquare$ {\parfillskip=0pt \finalhyphendemerits=0 \par}\medbreak}

\newcommand\capsize{\relax}
\newcommand\nobf{\noindent\bf}

\newcommand\diag{\operatorname{diag}}

\newcommand\interior{\operatorname{int}}

\newcommand\tr{\operatorname{tr}}

\newtheorem{lemma}{Lemma}[section]
\newtheorem{defi}[lemma]{Definition}
\newtheorem{theo}[lemma]{Theorem}
\newtheorem{prop}[lemma]{Proposition}
\newtheorem{coro}[lemma]{Corollary}

\title{The Asymptotics of Wilkinson's Shift: \\
Loss of Cubic Convergence}
\author{Ricardo S. Leite, Nicolau C. Saldanha and Carlos Tomei }

\begin{document}
\maketitle

\begin{abstract}
One of the most widely used methods for eigenvalue computation
is the $QR$ iteration with Wilkinson's shift:
here the shift $s$ is the eigenvalue of the bottom $2\times 2$
principal minor closest to the corner entry.
It has been a long-standing conjecture that the rate
of convergence of the algorithm is cubic.
In contrast, we show that there exist matrices for which
the rate of convergence is strictly quadratic.
More precisely, 
let $T_\xxxx$ be the $3 \times 3$ matrix having only two nonzero
entries $(T_\xxxx)_{12} = (T_\xxxx)_{21} = 1$
and let $\ILa$ be the set of real, symmetric
tridiagonal matrices with the same spectrum as $T_\xxxx$.
There exists a neighborhood $\bUu \subset \ILa$ of $T_\xxxx$
which is invariant under Wilkinson's shift strategy
with the following properties.
For $T_0 \in \bUu$, the sequence of iterates $(T_k)$ exhibits
either strictly quadratic or
strictly cubic convergence to zero of the entry $(T_k)_{23}$.
In fact, quadratic convergence occurs exactly when
$\lim T_k = T_\xxxx$.
Let $\bXx$ be the union of such quadratically convergent sequences $(T_k)$:
the set $\bXx$ has Hausdorff dimension $1$ and is
a union of disjoint arcs $\bXx^\sigma$ meeting at $T_\xxxx$,
where $\sigma$ ranges over a Cantor set.
\end{abstract}

\medbreak


\smallbreak

{\noindent\bf MSC-class:} 65F15; 37E99; 37N30.

\section{Introduction}

The $QR$ iteration
is a standard algorithm to compute eigenvalues of matrices in $\Tt$,
the vector space of $n \times n$ real symmetric tridiagonal matrices
(\cite{Wilkinson}, \cite{Demmel}, \cite{Parlett}).
More precisely, consider $T \in \Tt$  and a shift $s \in \RR$
so that $T - sI$ is an invertible matrix.
Given the $QR$ decomposition $T-sI = QR$,
where $Q$ is orthogonal and $R$ is upper triangular with positive diagonal,
the \textit{shifted step} obtains 
a new matrix $\bF(s,T) = Q^\transpose T Q = R T R^{-1}$.
Let $\bomega_-(T) \le \bomega_+(T)$
be the eigenvalues of the bottom $2 \times 2$ principal minor of $T$
and let $\bomega(T)$ be the eigenvalue closer to the bottom entry $(T)_{nn}$.
\textit{Wilkinson's shift} is the choice $s = \bomega(T)$ and
\textit{Wilkinson's step} is $\bW(T) = \bF(\bomega(T),T)$,
provided $\bomega(T)$ is well defined and is not an eigenvalue of $T$.

A matrix $T \in \Tt$ is \textit{unreduced} if $(T)_{i,i+1} = (T)_{i+1,i} \ne 0$
for $1 \le i < n$.
Recall that if $T_0$ is unreduced and the iterates $T_k = \bW^k(T_0)$, $k \in \NN$,
are well defined then $T_k$ is also unreduced and
the bottom off-diagonal entry $(T_k)_{n-1,n}$ tends to $0$.
The quick convergence of Wilkinson's algorithm
is a well known fact, as discussed in Section 8-11 of \cite{Parlett}:
we are interested in the precise rate of convergence of this sequence.
It has been conjectured (\cite{Parlett}, \cite{Demmel})
that it should be \textit{cubic}, i.e.,
$|(T_{k+1})_{n,n-1}| = O(|(T_{k})_{n,n-1}|^3)$.
Here, instead, we show that cubic convergence does not hold in general:
there are unreduced matrices $T_0$
for which the rate of convergence of the sequence
$(T_k)_{23}$ to $0$ is, in the words of Parlett, merely quadratic.

Given $T \in \Tt$, the above formulae imply that $\bW(T)$, if well defined,
is symmetric and upper Hessenberg and therefore
$T$ and $\bW(T)$ are matrices in $\Tt$ with the same spectrum.
For $\Lambda = \diag(1,-1,0)$,
denote by $\ILa$ the set of matrices in $\Tt$ similar to $\Lambda$
and consider $\bW$ as a map from $\ILa$ to $\ILa$.
There are some technical aspects to consider.
First, there are matrices $T \in \ILa$
with bottom entry $T_{3,3}$ equidistant from $\bomega_-$ and $\bomega_+$.
This introduces a step discontinuity in the map $\bW$:
when $T$ tends to such a matrix $T_0$,
$\bW(T)$ may approach either $\bF(\bomega_-(T_0),T_0)$
or $\bF(\bomega_+(T_0),T_0)$.
Also, strictly speaking, 
the definition of $\bW$ does not apply to matrices $T$ for which
$T_{23}=0$ because then $\bomega(T)$ is an eigenvalue of $T$.
It turns out, however, that $\bW$ can be continuously extended
to such matrices.

We introduce the notation required to state the main result of this paper.
Let \[ T_\xxxx =
\begin{pmatrix} 0 & 1 & 0 \\ 1 & 0 & 0 \\ 0 & 0 & 0 \end{pmatrix} \in \ILa. \]
An \textit{infinite sign sequence} is a function $\sigma: \NN \to \{+,-\}$.
Let $\Sigma$ be the set of infinite sign sequences:
$\Sigma$ admits the natural metric
$d(\sigma_1,\sigma_2) = 2^{-n}$ where $n$ is the smallest
number for which $\sigma_1(n) \ne \sigma_2(n)$.
The metric space $\Sigma$ is homeomorphic to the
middle-third Cantor set contained in $[0,1]$.
For $\sigma \in \Sigma$, let $\sigma^\sharp \in \Sigma$
be obtained by deleting the first sign, i.e.,
$\sigma^\sharp(n) = \sigma(n+1)$.

\begin{theo}
\label{theo:main}
There is an open neighborhood $\bUu \subset \ILa$ of $T_\xxxx$
with the following properties.
\begin{enumerate}[(a)]
\item{If $T \in \bUu$ then $\bW(T) \in \bUu$.}
\item{If $T_0 \in \bUu$ then the sequence $T_k = \bW^k(T_0)$
converges to $T_\infty$ with $(T_\infty)_{23} = 0$.
If $T_\infty = T_\xxxx$ then the convergence rate of $(T_k)_{23}$ to $0$
is strictly quadratic; otherwise it is strictly cubic.}
\item{Let $\bXx \subset \bUu$ be the set of initial conditions $T_0$
for which $T_\infty = T_\xxxx$.
Then $\bXx$ is the union of arcs $\bXx^\sigma$, $\sigma \in \Sigma$,
which are disjoint except for the common point $T_\xxxx$.
Also, $\bW(\bXx^\sigma) \subseteq \bXx^{\sigma^\sharp}$.}
\item{The set $\bXx \subset \bUu$ has Hausdorff dimension $1$.}
\end{enumerate}
\end{theo}


We opted to make this paper as self-contained as reasonably possible
at the expense of rendering several sections more technical.
The inductive arguments employed to control the iteration
make use of an assortment of explicit estimates.
The parameters obtained are loose and might be replaced
by a chain of existential arguments:
we hope that our choice led to a clearer presentation.

In Section 2 we introduce an explicit chart $\phi: \RR^2 \to \ILa$,
satisfying $\phi(p_\xxxx) = T_\xxxx$
where $p_\xxxx = (2,0)$.
Shifted steps in these $(x,y)$ coordinates, i.e., the functions
$F(s,x,y) = \phi^{-1}(\bF(s,\phi(x,y)))$,
admit a simple formula:
\[ F(s,x,y) = \left( \frac{1+s}{1-s} x, \frac{|s|}{1+s} y \right). \]
One may think of $(x,y)$ coordinates as a variation of spectral data
in the sense of \cite{BG}, \cite{DNT} and \cite{GH}:
eigenvalues and the absolute values of the first coordinates
of the normalized eigenvectors.
The map $\phi$ is an instance
of \textit{bidiagonal coordinates} (\cite{LST1}),
which form an atlas of $\ILa$:
$\phi$ is a chart whose image is an open dense set containing $T_\xxxx$.
In this chart, convergence issues reduce to local theory.
The present text might serve as an illustration of the general construction:
the downside is that the formula for $\phi$ is presented
without a natural process leading to it.

The rather technical 
Section 3 is dedicated to the study of the sub-eigenvalues
$\bomega_\pm(T)$ and their counterparts in $(x,y)$ coordinates
$\omega_\pm = \bomega_\pm \circ \phi$.
We also introduce the rectangle
$\Rr_a = [2-a,2+a]\times [-a/10,a/10]$, $0 < a \le 1/10$,
whose image $\bRr_a = \phi(\Rr_a)$ is the closure of the invariant
neighborhood $\bUu$ in Theorem \ref{theo:main}.
The discontinuities of $\omega$ lie on the vertical line $x = 2$:
if $x < 2$ (resp. $x > 2$), $\omega(x,y) = \omega_+(x,y)$
(resp. $\omega(x,y) = \omega_-(x,y)$).

In Section 4 we show that $\Rr_a$ is invariant under
$W = \phi^{-1} \circ \bW \circ \phi$ (for sufficiently small $a$).
The discontinuous map $W$ has smooth restrictions
to the interior of the rectangles
$\Rr_{a,+} = [2-a,2] \times [-a/10,a/10]$ and
$\Rr_{a,-} = [2,2+a] \times [-a/10,a/10]$
which extend continuously to
$W_\pm: \Rr_{a,\pm} \to \Rr_a$.

The rest of the proof of Theorem \ref{theo:main}
proceeds by studying the iterations of $W$.
At this point, the vocabulary and techniques
from dynamical systems are natural
(another application of dynamical systems to numerical spectral theory
is the work of Batterson and Smillie \cite{BS}
on the Rayleigh quotient iteration).
In Section 5 we provide different characterizations
of the set $\Xx = \phi^{-1}(\bXx)$ and show its non-triviality.

In Section 6 we construct a continuous map from
$[-a/10,a/10] \times \Sigma \to \Xx$
taking a pair $(y,\sigma)$ to the only point $p_0 = (g^\sigma(y),y) \in \Xx$
for which $p_k = W^k(p_0) \in \Rr_{a,\sigma(k)}$.
Informally, the sign sequence $\sigma$ specifies the side of the rectangle
$\Rr_a$ in which $p_k$ lies.
This yields the decomposition of $\Xx$ as
a union of the Lipschitz arcs $\Xx^\sigma$, $\sigma \in \Sigma$.
Sharper estimates concerning $g^\sigma(y)$ are needed to prove that
$\Xx$ is very thin, i.e., that its Hausdorff dimension is $1$.
The concept of Hausdorff dimension is only used at the very end
of the argument in order to obtain a more concise formulation
of an estimate on the number of balls of radius $r$ required to cover $\Xx$.

\section{Local coordinates and $s$-steps}

The {\it $QR$ decomposition}
of an invertible matrix $M$ is $M = \bQ(M)\,\bR(M)$,
where $\bQ(M)$ is orthogonal and $\bR(M)$ is upper
triangular with nonnegative diagonal.
The \textit{$s$-step} with shift $s$ is the map
\begin{align}
\bF(s,M) &= (\bQ(M- sI))^\ast\,M\,\bQ(M - sI) \label{eq:stepQ} \\
&= \bR(M - sI)\,M\,(\bR(M - sI))^{-1}. \label{eq:stepR}
\end{align}
This is well defined provided $M-sI$ is invertible.

As an example, to be used in the sequel, consider
\[ T_\xxxx = \begin{pmatrix} 0&1&0\\1&0&0\\0&0&0 \end{pmatrix}. \]
Clearly, the matrix $T_\xxxx - sI$ is invertible if $s \ne 0, \pm 1$.
For $s \ne 0$, $|s| < 1$,
the $QR$ decomposition of $T_\xxxx - sI$ is
\[ \bQ(T_\xxxx - sI)\ \bR(T_\xxxx - sI) = \begin{pmatrix}
{-s s_1} & {s_1} & {0} \\
{s_1} & {s s_1} & {0} \\
{0} & {0} & {-\sign(s)}
\end{pmatrix}
\begin{pmatrix}
{s_1^{-1}} & {-2s s_1} & {0} \\
{0} & {(1-s^2)s_1} & {0} \\
{0} & {0} & {|s|}
\end{pmatrix} \]
where $s_1 = (1+s^2)^{-1/2}$.
Notice that $\bQ(T-sI)$ has a jump discontinuity at $s=0$
while $\bR(T_\xxxx - sI)$ is continuously defined but
ceases to be invertible at $s = 0$.
Standard $s$-steps for $T_\xxxx$ are given by
\[  \bF(s,T_\xxxx) = \frac{1}{1+s^2} \begin{pmatrix}
{-2s} & {1-s^2} & {0} \\
{1-s^2} & {2s} & {0} \\
{0} & {0} & {0}
\end{pmatrix}. \]
Notice that this formula extends continuously to $s=0$
with $\bF(0,T_\xxxx) = T_\xxxx$.
Such continuous extensions will be important throughout the paper.


Let $\Tt$ be the set of $3 \times 3$ real, symmetric, tridiagonal matrices.
Let $\Lambda= \diag(1,-1,0)$
and $\ILa = \{ Q^\ast\Lambda Q, Q \in O(n) \}$ be
the set of matrices in $\Tt$ similar to $\Lambda$.
In fact, $\ILa$ is a bitorus (\cite{Tomei}, \cite{LST1})
but this kind of global information will not be used in this paper.

We now define the relevant parametrization of $\ILa$ near $T_\xxxx$.
Let $\Ee$ be the set of \textit{signature matrices}, i.e.,
diagonal $n \times n$ matrices with diagonal entries equal to $\pm 1$.
An invertible matrix $M$ is \textit{$LU$-positive}
if there exist lower and upper triangular matrices $L$ and $U$
with positive diagonals so that $M = LU$;
equivalently, $M$ is $LU$-positive if its leading principal minors
have positive determinant.
For instance, the only $LU$-positive orthogonal matrix $Q$ with
$T_\xxxx = Q^\ast\Lambda Q$ is
\[ Q = Q_\xxxx = \begin{pmatrix} 1/\sqrt2 & 1/\sqrt2 & 0
\\ -1/\sqrt2 & 1/\sqrt2 & 0 \\
0 & 0 & 1 \end{pmatrix}.  \]

\begin{prop}
\label{prop:phi}
Set $p_\xxxx = (2,0)$ and $\phi: \RR^2 \to \Tt$,
\[ \phi(x,y) =  \frac{1}{r_1^2r_2^2}
\begin{pmatrix}
(4 - x^2) r_2^2 & 2x r_2^3 & 0 \\
2x r_2^3 & -4(4 - x^2 - 4x^2y^4 + x^4y^4) & 2y r_1^3 \\
0 & 2y r_1^3 & y^2(x^2 - 4) r_1^2
\end{pmatrix}, \]
where $r_1 = \sqrt{4 + x^2 + 4x^2y^2}$ and $r_2 = \sqrt{4 + 4y^2 + x^2y^2}$.
Then $\phi(p_\xxxx) = T_\xxxx$ and
there exist neighborhoods $\Uu_1 \subset \RR^2$ of $p_\xxxx$
and $\Bb_1 \subset \Tt$ of $T_\xxxx$
such that $\phi|_{\Uu_1}$ is an immersion
with image $\bUu_1 = \ILa \cap \Bb_1$.
Also, $\sign(x) = \sign((\phi(x,y))_{12})$
and $\sign(y) = \sign((\phi(x,y))_{23})$.
Finally, for matrices $T = \phi(x,y) \in \bUu_1$,
the only $LU$-positive orthogonal matrix
for which $T = Q^\ast \Lambda Q$ is
\[ Q =
\frac{1}{r_1r_2}
\begin{pmatrix}
{2}{r_2} & {2x(1+2y^2)} & {xy}{r_1} \\
{-x}{r_2} & {2(2+x^2y^2)} & {-2y}{r_1} \\
{-2xy}{r_2} & {y(4-x^2)} & {2}{r_1}
\end{pmatrix}. \]
\end{prop}

Since we are interested in the behavior of Wilkinson's step
near $T_\xxxx$, we successively define nested neighborhoods
$\bUu_k$ of $T_\xxxx$ and $\Uu_k$ of $p_\xxxx$:
the process stops with the contruction
of the rectangle $\Rr_a$ at Proposition \ref{prop:split}.
Constructions on matrices use boldface symbols.

\proof
The example above obtains $T_\xxxx = \phi(p_\xxxx)$.
The remaining statements also follow from numerical checks, but instead
we provide  motivation for the construction of the map $\phi$.
Given the eigenvalue matrix $\Lambda$, any matrix $T \in \ILa$ admits an
orthogonal diagonalization $T = Q^\ast  \Lambda  Q$, which is not unique:
all other orthogonal diagonalizations are of the form
$T = Q^\ast E \Lambda EQ$ for some signature matrix $E \in \Ee$.

Now, the rows of $Q$ are normal eigenvectors of $T$,
so its columns $q_1, q_2$ and $q_3$ are also orthonormal.
For vectors $q$ and $w$, we write $q \sim w$ to indicate
that they are collinear.
Say $q_1 \sim w_1 = (2, -x,-2xy)^\ast$.
It is well known that simple eigenvalues and their normalized eigenvectors
vary smoothly with the related matrix;
thus, for $T \in \ILa$, $T$ near $T_\xxxx$,
the first column of $Q$ can indeed be written in the form above
for an appropriate choice of $x$ and $y$,
$(x,y)$ near $p_\xxxx$.
Since $0= T_{31} = \langle q_3, \Lambda q_1  \rangle$, we
must have $q_3$ orthogonal both to $w_1$ and $\Lambda w_1$
and therefore $q_3 \sim w_1 \times \Lambda w_1 \sim w_3 = (xy, -2y, 2)^\ast$.
Similarly, $q_2 \sim w_2 = w_3 \times w_1 =
(2x(1+2y^2),2(2+x^2y^2),y(4-x^2))^\ast$.
Thus, we search for an $LU$-positive orthogonal matrix $Q = EWN$
where $E \in \Ee$, $W$ has columns $w_i$
and $N$ is a positive diagonal matrix.
It is now easy to verify that $N = \diag(1/r_1,1/(r_1r_2),1/r_2)$
where $r_1$ and $r_2$ are given in the statement;
also, by computing signs of pricipal minors,
$WN$ is $LU$-positive and therefore $E = I$.
Expansion of the product $Q^\ast\Lambda Q$ obtains the formula for $\phi$.

To show that $\phi$ is an immersion near $T_\xxxx$,
it suffices to verify that the Jacobian $D \phi$ at $p_\xxxx$ is injective.
This is easily seen by
computing partial derivatives of the entries $T_{11}$ and $T_{23}$
of $T = \phi(x,y)$ with respect to $x$ and $y$.
\qed

Actually, $\phi$ is a diffeomorphism from $\RR^2$ to its image \cite{LST1},
but this will not be used in this text.
A reason for choosing such a parametrization for the first column of $Q$
will be clarified by the next proposition:
an $s$-step admits a simple representation in $(x,y)$-coordinates.

\begin{prop}
\label{prop:s-step}
There are open neighborhoods $S_2 \subset (-1,1)$ of $s = 0$,
$\Uu_2 \subset \Uu_1$ of $p = p_\xxxx$ and
$\bUu_2 = \phi(\Uu_2) \subset \bUu_1$ with the following properties.
\begin{enumerate}[(a)]
\item{The function $\bF: (S_2 \smallsetminus \{0\}) \times \bUu_2 \to \bUu_1$
 is smooth and
extends continuously (but not smoothly) to $S_2 \times \bUu_2$.}
\item{Let $F: S_2 \times \Uu_2 \to \Uu_1$
be $\bF$ expressed in $(x,y)$-coordinates, i.e.,
$F(s,p) = \phi^{-1}(\bF(s,\phi(p)))$.
Then, for $p = (x,y) \in \Uu_2$,
\[ F(s,x,y) = \left( \frac{1+s}{1-s} x, \frac{|s|}{1+s} y \right). \]}
\end{enumerate}
\end{prop}

\proof
By definition, $\bF(s,T) = (\bQ(T - sI))^\ast  T  \bQ(T - sI)$.
For $T \approx T_\xxxx$ (i.e., $T$ near $T_\xxxx$)
and $s \approx 0$, $s \ne 0$, we have $T - sI \approx T_\xxxx$.
The first two columns of $\bQ(T - sI)$ can be obtained from
the corresponding columns of $T - sI$ by the Gram-Schmidt algorithm
and therefore
\[ \bQ(T - sI) \approx
\begin{pmatrix} 0 & 1 & 0 \\ 1 & 0 & 0 \\
0 & 0 & -\sign(\det(T-sI)) \end{pmatrix}. \]
Whatever the sign of $\det(T-sI)$, $\bF(s,T) \approx T_\xxxx$ so that
\[ \lim_{s \to 0, s \ne 0, T \to T_\xxxx} \bF(s,T) = T_\xxxx. \]
This proves that there exist open neighborhoods $S_2$ and $\bUu_2$
such that $\bF: (S_2 \smallsetminus \{0\}) \times \bUu_2 \to \bUu_1$
is well defined and extends continuously to the point $(0,T_\xxxx)$
with $\bF(0,T_\xxxx) = T_\xxxx$.
The smoothness of $\bF$ in $(S_2 \smallsetminus \{0\}) \times \bUu_2$
follows from the smoothness of $\bQ$ in the set of invertible matrices.
The fact that $\bF$ extends continuously to $S_2 \times \bUu_2$
will follow from the explicit formula in item (b).

As in the proof of Proposition \ref{prop:phi},
for $(s,T_0) \in S_2 \times \bUu_2$,
write $T_0 - sI = Q_0^\ast (\Lambda - sI) Q_0$,
where $Q_0$ is $LU$-positive.
Notice that if $(a_0,b_0,c_0)^\ast$ is the first column of $Q_0$
and $T_0 = \phi(x_0,y_0)$ then $x_0 = -2b_0/a_0$ and $y_0 = c_0/(2b_0)$.
For $s \ne 0$, let $(x_1,y_1) = G(s,x_0,y_0)$
and $T_1 = \bF(s,T_0) = \phi(x_1,y_1) = Q_1^\ast \Lambda Q_1$
where $Q_1$ is $LU$-positive.
Let $(a_1,b_1,c_1)$ be the first column of $Q_1$:
we have $x_1 = -2b_1/a_1$ and $y_1 = c_1/(2b_1)$.
We must therefore compute $a_1, b_1, c_1$.

By definition,
\[ T_1 = (\bQ(T_0 - sI))^\ast  T_0  \bQ(T_0 - sI)
= (\bQ(T_0 - sI))^\ast  Q_0^\ast  \Lambda  Q_0 \bQ(T_0 - sI),\]
so that $Q_1 = E Q_0 \bQ(T_0 - sI)$
for some signature matrix $E = \diag(e_1, e_2, e_3) \in \Ee$.
Set $Q_1 = E \bQ(Q_0(T_0 - sI)) = E \bQ((\Lambda - sI) Q_0)$.
Assuming $|s| < 1$, we have the positive collinearity
$(a_1,b_1,c_1)^\ast \sim (e_1 (1-s) a_0, e_2 (-1-s) b_0, e_3 (-s) c_0)^\ast$
and thus
\[ x_1 = -\frac{e_2}{e_1} \frac{1+s}{1-s} x_0, \quad
y_1 = \frac{e_3}{e_2} \frac{s}{1+s} y_0. \]
Now, $\sign(x_i) = \sign((T_i)_{12})$, $i = 0, 1$,
and $\sign((T_1)_{12}) = \sign((T_0)_{12})$
(from Proposition \ref{prop:phi} and equation (\ref{eq:stepR}))
and therefore $\sign(x_1) = \sign(x_0)$;
similarly $\sign(y_1) = \sign(y_0)$, completing the proof.
\qed

The computations above fit into a larger context, which we now outline.
As with $s$-steps, many eigenvalue algorithms act on $n \times n$
\textit{Jacobi matrices} (as in \cite{BG}, \cite{GH} and \cite{Moser},
real symmetric tridiagonal matrices $T$
with $T_{i+1,i} > 0$, $i = 1, 2, \ldots, n-1$)
with non-Jacobi limit points.
Recall that Jacobi matrices have simple (real) eigenvalues $\lambda_i$
and that its normalized eigenvectors $v_i$ can be chosen so that
the first coordinates $c_i$ are positive; notice that $\sum c_i^2 = 1$.
The eigenvalues $\lambda_i$ and the \textit{norming constants} $c_i$
form a standard set of coordinates for Jacobi matrices.
In these standard coordinates (\cite{DNT}),
the $s$-step acting on a Jacobi matrix $T$ is given by
\[(\lambda_i, c_i), \quad \mapsto \quad
(\lambda_i,\frac{ |\lambda_i-s| c_i}{\sum_i (|\lambda_i-s| c_i)^2}),
\quad i=1,\ldots,n \]
provided $T - sI$ is invertible.
Such coordinates require some modification
in order to extend beyond the set of Jacobi matrices.
The {\it bidiagonal coordinates} in \cite{LST1} are,
up to multiplicative constants, quotients $c_{i+1}/c_i$
which admit a simpler evolution under $s$-steps.
The map $\phi(x,y)$ retrieves a matrix from
one among $3!$ possible choices of bidiagonal coordinates.
In general, each permutation $\pi$ of the eigenvalues
obtains a map $\phi_\pi: \RR^n \to \ILa$ and this family of maps
is an atlas for $\ILa$.


\section{Sub-eigenvalues}

For a matrix $T \in \Tt$,
let $\hat T$ be its bottom $2 \times 2$ diagonal block.
Given $T \in \ILa$,
the eigenvalues $\bomega_-(T) \le \bomega_+(T)$ of $\hat T$
are the {\em sub-eigenvalues} of $T$.
Notice that $\bomega_-(T_\xxxx) = \bomega_+(T_\xxxx) = 0$.

Consider the circles $\bc_h, \bc_v \subset \ILa$ through $T_\xxxx$
parametrized by
\[ \begin{pmatrix}
\sin\theta & \cos\theta & 0 \\
\cos\theta & -\sin\theta & 0 \\
0 & 0 & 0
\end{pmatrix}, \quad
\begin{pmatrix}
{0} & {\cos\theta} & {0} \\
{\cos\theta} & {0} & {\sin\theta} \\
{0} & {\sin\theta} & {0}
\end{pmatrix}, \]
respectively.
Clearly, $T \in \bc_h$ if and only if $T_{23} = T_{33} = 0$;
also, $T \in \bc_v$ if and only if $T_{11} = T_{22} = T_{33} = 0$.

\begin{prop}
\label{prop:omega}
There is an open connected neighborhood $\bUu_3 \subset \bUu_2$ of $T_\xxxx$
with the following properties.
\begin{enumerate}[(a)]
\item For all $T \in \bUu_3$ the following interlacing inequalities hold:
\[ -1 \le \bomega_-(T) \le 0 \le \bomega_+(T) \le 1, \quad
\bomega_-(T) \le T_{33} \le \bomega_+(T). \]
\item The only matrix $T \in \bUu_3$ for which $\bomega_-(T) = \bomega_+(T)$
is $T = T_\xxxx$.
The functions $\bomega_\pm: \bUu_3 \to \RR$ are continuous
and smooth on $\bUu_3 \smallsetminus \{T_\xxxx\}$.
\item
A matrix $T \in \bUu_3$ belongs to $\bc_h$
if and only if at least one sub-eigenvalue of $T$
coincides with an eigenvalue of $T$.
In particular, $T_{12} > 0$ for all $T \in \bUu_3$.
\item
A matrix $T \in \bUu_3$ belongs to $\bc_v$
if and only if $T_{33}$ is equidistant
from the sub-eigenvalues of $T$.
\item
For all $T \in \bUu_3$, $\bomega_\pm(T) \in S_2$.
\end{enumerate}
\end{prop}

The set $S_2$ mentioned in item (e)
was defined in Proposition \ref{prop:s-step}.

\proof
The first inequality is the interlacing of the eigenvalues of $T$ and $\hat T$,
the second is the interlacing of those of $\hat T$ and $T_{33}$.

From item (a), if the sub-eigenvalues are equal then
$\bomega_-(T) = \bomega_+(T) = T_{33} = 0$.
Thus, $\hat T = 0$ and since the trace of $T$ equals $0$,
one has $T_{11}=0$ and $T_{12}= \pm 1$:
only the positive choice, which gives rise to $T_\xxxx$ itself, is relevant.
Continuity in $\bUu_3$ and smoothness in $\bUu_3 \smallsetminus \{T_\xxxx\}$
of the functions $\bomega_\pm$ are now easy.

For item (c), suppose that $\det(T - \bomega_+(T) I)=0$
(the case $\det(T - \bomega_-(T) I) = 0$ is similar).
By construction, $\bomega_+(T)$ must be a common
eigenvalue of $T$ and $\hat T$.
Expand the characteristic polynomial of $T$
along the first row to obtain
\[ \det(\lambda I - T) =
(\lambda - T_{11}) \det(\lambda I - {\hat T}) +
( - T_{12}^2 \lambda + T_{12}^2 T_{33}). \]
A common eigenvalue annihilates the two determinants,
thus $T_{12}^2(T_{33} - \bomega_+(T))=0$.
Since $T_{12} \approx 1$ and $\bomega_+(T)$ equals an eigenvalue of $\Lambda$
we have $\bomega_+(T) = T_{33} = 0$.
Now $\det\hat T = 0$, which implies $T_{23}=0$.
Notice that $T_{12} = 0$, $T_{23} \ne 0$ implies that
some sub-eigenvalue equals $\pm 1$:
this possibility was excluded above.

For item (d),
if $T_{33}$ is equidistant from the sub-eigenvalues,
$T \in \ILa$ must be of the form
\[ \begin{pmatrix}
{-2d} & b & {0} \\
b & d & c\\
{0} & c & d
\end{pmatrix}.\]
Now $\det(T - \lambda)= \lambda^3 - \lambda
= \lambda^3+(-3d^2-b^2-c^2)\lambda-2dc^2+db^2+2d^3$, so
$d(b^2 - 2 c^2 + 2 d^2) = 0$ and $b^2 + c^2 + 3 d^2 = 1$.
If $d=0$, $T$ belongs to $\bc_v$:
notice that $T_\xxxx$ corresponds to $\theta = 0$
in the parametrization of $\bc_v$.
If $d \ne 0$, $T$ lies in one of two curves in $\ILa$
which may be assumed to be disjoint from $\bUu_3$.

Continuity of the functions $\bomega_\pm$
allows for a choice of $\bUu_3$ satisfying the final condition.
\qed

Let $\Uu_3 = \phi^{-1}(\bUu_3)$.
From item (c) above, $(x,y) \in \Uu_3$ implies $x > 0$.
Let $r_h \subset \RR^2$ be the horizontal axis and
$r_v \subset \RR^2$ be the vertical line $x = 2$.
Let
$\Uu_{3,+} = \Uu_3 \cap \{(x,y) \; | \; x \le 2 \}$,
$\Uu_{3,-} = \Uu_3 \cap \{(x,y) \; | \; x \ge 2 \}$,
$\bUu_{3,\pm} = \phi(\Uu_{3,\pm})$.

It is convenient to work with a domain for $(x,y)$ coordinates
which is more explicit than the set $\Uu_3$.

\begin{defi}
Let $\Rr_a = [2-a,2+a] \times [-\ver, \ver] \subset \Uu_3$
be a rectangle centered in $p_\xxxx=(2,0)$.
The rectangle $\Rr_a$ is split by $r_v$ in two closed rectangles
$\Rr_{a,+} = \Rr_a \cap \Uu_{3,+}$ and $\Rr_{a,-} = \Rr_a \cap \Uu_{3,-}$.
\end{defi}

From \cite{LST1}, it follows easily that $\Uu_3$
can be taken to contain $\Rr_{1/10}$:
using this result, as we shall see, all the subsequent constructions
are compatible with $a = 1/10$.
To make this paper self-contained the working hypothesis
is only $a \le 1/10$, $\Rr_a \subset \Uu_3$.

We rephrase Proposition \ref{prop:omega} in $(x,y)$-variables.
Write $\omega_\pm(x,y) = \bomega_\pm(\phi(x,y))$
so that the functions $\omega_\pm$ are continuous in $\Rr_a$
with a non-smooth point $p_\xxxx$.
Set $\bRr_a = \phi(\Rr_a)$.

\begin{prop}
\label{prop:split}
The diffeomorphism $\phi: \Rr_a \to \bRr_a \subset \bUu_3$
yields bijections from $r_h \cap \Rr_a$ to $\bc_h \cap \bRr_a$ and
from $r_v \cap \Rr_a$ to $\bc_v \cap \bRr_a$.
The functions $\omega_\pm(x,y)$ are even with respect to $y$, i.e,
$\omega_\pm(x,y) = \omega_\pm(x,-y)$.
For points $(x,y) \in \Rr_{a,+}$ (resp. $\Rr_{a,-}$),
$(\phi(x,y))_{33}$ is to the right (resp. left) of
$(\omega_+(x,y) + \omega_-(x,y))/2$.
\end{prop}

Signs in the notation $\Rr_{a,\pm}$
indicate which among $\omega_\pm$ is closer to $(\phi(x,y))_{33}$:
unfortunately, they are the reverse of what their position might suggest.

\proof
We already saw in Proposition \ref{prop:phi} that
$\sign((\phi(x,y))_{23}) = \sign(y)$.
Clearly, $\phi(2,(\tan\theta)/\sqrt{2})$ equals the matrix used
to parametrize $\bc_v$ in the statement of Proposition \ref{prop:omega}.
Evenness of $\omega_\pm(x,y)$ in $y$
is immediate from the explicit form of $\phi$ in Proposition \ref{prop:phi}.
Finally, to decide which sub-eigenvalue of $T = \phi(x,y)$
is closer to $T_{33}$,
compare $\tr \hat T = T_{22} + T_{33} =
\omega_+(x,y) + \omega_-(x,y)$ with $2T_{33}$:
\[ (2T_{33} - \tr T) r_1^2 r_2^2 = (x-2)(x+2)(x^2y^2+8x^2y^4-4+4y^2). \]
For $x$ slightly smaller (resp. larger) than $2$,
the expression is positive (resp. negative) and thus
$T_{33}$ is to the right (resp. left) of
$(\omega_+(x,y) + \omega_-(x,y))/2$.
Since the only points where $T_{33} = (\omega_+(x,y) + \omega_-(x,y))/2$
are those in $r_v \cap \Rr_a$ and $\Rr_a$ is connected the result follows.
\qed

We need more precise estimates for the sub-eigenvalues $\omega_\pm$
near $p_\xxxx = (2,0)$.

\begin{defi}
The \textit{wedge} $\Xx_0$ is
$\{(x,y) \in \Rr_a \;|\; |y| \ge |x-2|/10 \}$;
set $\Xx_{0,\pm} = \Xx_0 \cap \Rr_{a,\pm}$.
\end{defi}

\begin{figure}[ht]
\psfrag{2-a}{$2-a$}
\psfrag{2+a}{$2+a$}
\psfrag{Ra+}{$\Rr_{a,+}$}
\psfrag{Ra-}{$\Rr_{a,-}$}
\psfrag{D0+}{$\Dd_{0,+}$}
\psfrag{D0-}{$\Dd_{0,-}$}
\psfrag{X0+}{$\Xx_{0,+}$}
\psfrag{X0-}{$\Xx_{0,-}$}
\psfrag{rh(y=0)}{$r_h$ ($y=0$)}
\psfrag{rv(x=2)}{$r_v$ ($x=2$)}
\psfrag{pd}{$p_\xxxx$}
\begin{center}
\epsfig{height=35mm,file=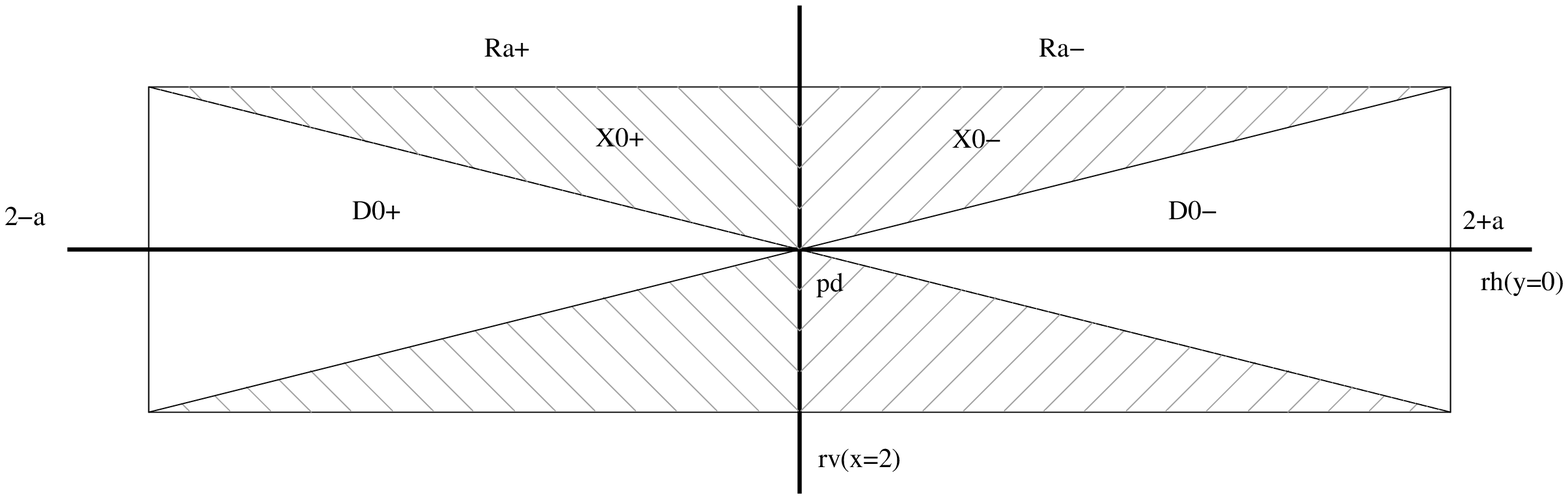}
\end{center}
\caption{The rectangle $\Rr_a$;
in $\Rr_{a,+}$, $\omega = \omega_+$ and
in $\Rr_{a,-}$, $\omega = \omega_-$.}
\label{fig:tab0}
\end{figure}

Figure \ref{fig:tab0} contains some of the geometric objects
defined in this section; the triangles $\Dd_{0,\pm}$
will be defined in the next section.

\begin{lemma}
\label{lemma:omega}
Near the point $p_\xxxx = (2,0)$ the functions $\omega_\pm$
have a cone-like behavior:
\[ \omega_\pm = \frac{(x-2) \pm \sqrt{(x-2)^2 + 32y^2}}{4} + O(d^2) \]
where $d^2 = (x-2)^2 + y^2$.
In $\Rr_{a,\pm}$, $y^2 \le |\omega_\pm(x,y)| \le 2|y|$;
in $\Xx_{0,\pm}$, $|y|/5 \le |\omega_\pm(x,y)| \le 2|y|$.
\end{lemma}

\proof
As in Proposition \ref{prop:phi}, set $r_1^2 = 4 + x^2 + 4x^2y^2$.
The sub-eigenvalues solve
$r_1^2 \omega^2+(4-x^2)\omega-4x^2y^2=0$,
\begin{equation}
\omega_\pm = \frac{-4 + x^2 \pm \sqrt\Delta}{2r_1^2},
\label{eq:omega}
\end{equation}
where $\Delta = \Delta(x,y) = ((x+2)^2 + 8x^2y^2)((x-2)^2 + 8x^2y^2) \ge 0$,
$\Delta = O(d^2)$.
In particular, $\Delta = 0$ in $\Rr_a$ only for $p_\xxxx=(2,0)$.
The expression for $\omega_\pm$ in the statement of the lemma
follows directly from
\[ \Delta - 16((x - 2)^2 + 32y^2) = O(d^3).\]
For $y > 0$, the inequality $\omega_+ \le 2|y|$ is equivalent to
$\Delta \le \left(4r_1^2 y + 4 - x^2\right)^2$ which is in turn equivalent to
\[ 8 r_1^2 y (4 - x^2 + 8y + 8x^2y^3) \ge 0, \]
which is clearly true for $0 < x \le 2$.
The case $y < 0$ follows from the fact that $\omega_+$ is even in $y$.
A similar factorization holds for $\omega_-$.
Also, we have $\omega_+(x,y) \omega_-(x,y) = -4x^2y^2/r_1^2$
and $|\omega_-(x,y)|, |\omega_+(x,y)| \le 1$.
Since $a \le 1/10$ we have that $4x^2/r_1^2 \le 1$ for $(x,y) \in \Rr_a$.
The estimate for $\omega_\pm$ in $\Xx_{0,\pm}$ is now
somewhat cumbersome but straightforward.
\qed

We also need estimates for partial derivatives of $\omega_\pm$
in terms of $x$ and $y$.

\begin{lemma}
\label{lemma:partialomega}
For all $(x,y) \in \Rr_a - \{p_\xxxx\}$
the partial derivatives $(\omega_\pm)_x$ and $(\omega_\pm)_y$
satisfy $0 \le (\omega_\pm)_x < 1$ and $|(\omega_\pm)_y| < 7/3$.
The equality $(\omega_+)_x = 0$ (resp. $(\omega_-)_x = 0$)
holds exactly for $y = 0$, $x < 2$ (resp. $x > 2$);
for $y \ne 0$, $\pm y (\omega_\pm)_y > 0$.
\end{lemma}

\proof
The partial derivatives of $\omega_\pm$ are
\begin{align}
(\omega_\pm)_x &=
\frac{8x}{r_1^4 \sqrt{\Delta}}
\left( \left( (1+2y^2) \sqrt{\Delta} \right) \pm
\left(-4 + x^2 + 8y^2 + 6x^2y^2 + 16x^2y^4 \right) \right), \notag\\
(\omega_\pm)_y &=
\frac{4x^2y}{r_1^4 \sqrt{\Delta}}
\left( \left( (4-x^2) \sqrt{\Delta} \right) \pm
\left(16 + 24x^2 + x^4 + 32x^2y^2 + 8x^4y^2 \right) \right), \notag
\end{align}
which are well defined provided $\Delta \ne 0$,
i.e., outside of $p_\xxxx=(2,0)$.
Also,
\[ \left( (1+2y^2) \sqrt{\Delta} \right)^2 -
\left(-4 + x^2 + 8y^2 + 6x^2y^2 + 16x^2y^4 \right)^2
= 8 y^2 r_1^4 \ge 0 \]
whence
$(1+2y^2) \sqrt{\Delta} \ge
\left|-4 + x^2 + 8y^2 + 6x^2y^2 + 16x^2y^4 \right|$
and $(\omega_\pm)_x \ge 0$; equality implies $y = 0$.
Also, $(\omega_\pm)_x \le 16x(1+2y^2)/r_1^4 < 1$.


In the rectangle $\Rr_a$,
\[115 < 16 + 24x^2 + x^4 + 32x^2y^2 + 8x^4y^2 < 142, \quad
\left| (4-x^2) \sqrt{\Delta} \right| < 1/4 \]
and the signs of $(\omega_\pm)_y$ are settled.
Since
\[ \frac{y^2}{\Delta} \le
\frac{1}{(x+2)^2 + 8x^2y^2} \; \frac{y^2}{8 x^2 y^2} < \frac{1}{360} \]
we also have
\[ |(\omega_\pm)_y| \le \frac{|y|}{\sqrt{\Delta}} \;
\frac{4x^2 \cdot \left(
16 + 24x^2 + x^4 + 32x^2y^2 + 8x^4y^2 + \left| (4-x^2) \sqrt{\Delta} \right|
\right) }{r_1^4}
< \frac{7}{3}. \]
\qed



\section{Wilkinson's step}

Take $\bomega(T)$ to be the sub-eigenvalue of $T$ closer to $T_{33}$;
in case of a draw, we arbitrarily choose $\bomega(T) = \bomega_+(T)$.
{\em Wilkinson's step}
is the map $\bW(T) = \bF(\bomega(T),T)$.
From now on we shall work in $(x,y)$ coordinates, i.e.,
with $W(x,y) = \phi^{-1}(\bW(\phi(x,y))) = F(\omega(x,y),x,y)$
where $\omega(x,y) = \bomega(\phi(x,y))$.
Write $W_\pm(x,y) = F(\omega_\pm(x,y),x,y)$.
From Propositions \ref{prop:s-step} and \ref{prop:omega},
the maps $W_\pm: \Rr_a \to \Uu_1$
are continuous and $W: \Rr_a \to \Uu_1$ is well defined
with step discontinuities along $r_v \cap \Rr_a$.
From Proposition \ref{prop:s-step}
and the fact that $\omega_-(x,y) \le 0 \le \omega_+(x,y)$,
\begin{equation}
W_\pm(x,y) = (X_\pm(x,y),Y_\pm(x,y)) =
\left( \frac{1+ \omega_\pm(x,y)}{1-\omega_\pm(x,y)}\;x, 
\frac{\pm\omega_\pm (x,y)}{1 + \omega_\pm (x,y)}\;y \right)
\label{eq:W}
\end{equation}
and therefore $W_\pm$ are smooth functions in
$\Rr_a \smallsetminus \{p_\xxxx\}$.
Evenness of $\omega_\pm$ with respect to $y$
(Proposition \ref{prop:split})
yields $X_\pm(x,y) = X_\pm(x,-y)$, $Y_\pm(x,y) = -Y_\pm(x,-y)$.

\begin{figure}[ht]
\psfrag{-2e-4}{$\scriptscriptstyle -2 \cdot 10^{-4}$}
\psfrag{2e-4}{$\scriptscriptstyle 2 \cdot 10^{-4}$}
\psfrag{2.1}{$\scriptscriptstyle 2.1$}
\psfrag{2.04}{$\scriptscriptstyle 2.04$}
\psfrag{2}{$\scriptscriptstyle 2$}
\psfrag{1.96}{$\scriptscriptstyle 1.96$}
\psfrag{1.9}{$\scriptscriptstyle 1.9$}
\psfrag{Wp(Rp)}{$W_+(\Rr_{a,+})$}
\psfrag{Wm(Rm)}{$W_-(\Rr_{a,-})$}
\psfrag{W(R)}{$W(\Rr_{a})$}
\begin{center}
\epsfig{height=60mm,file=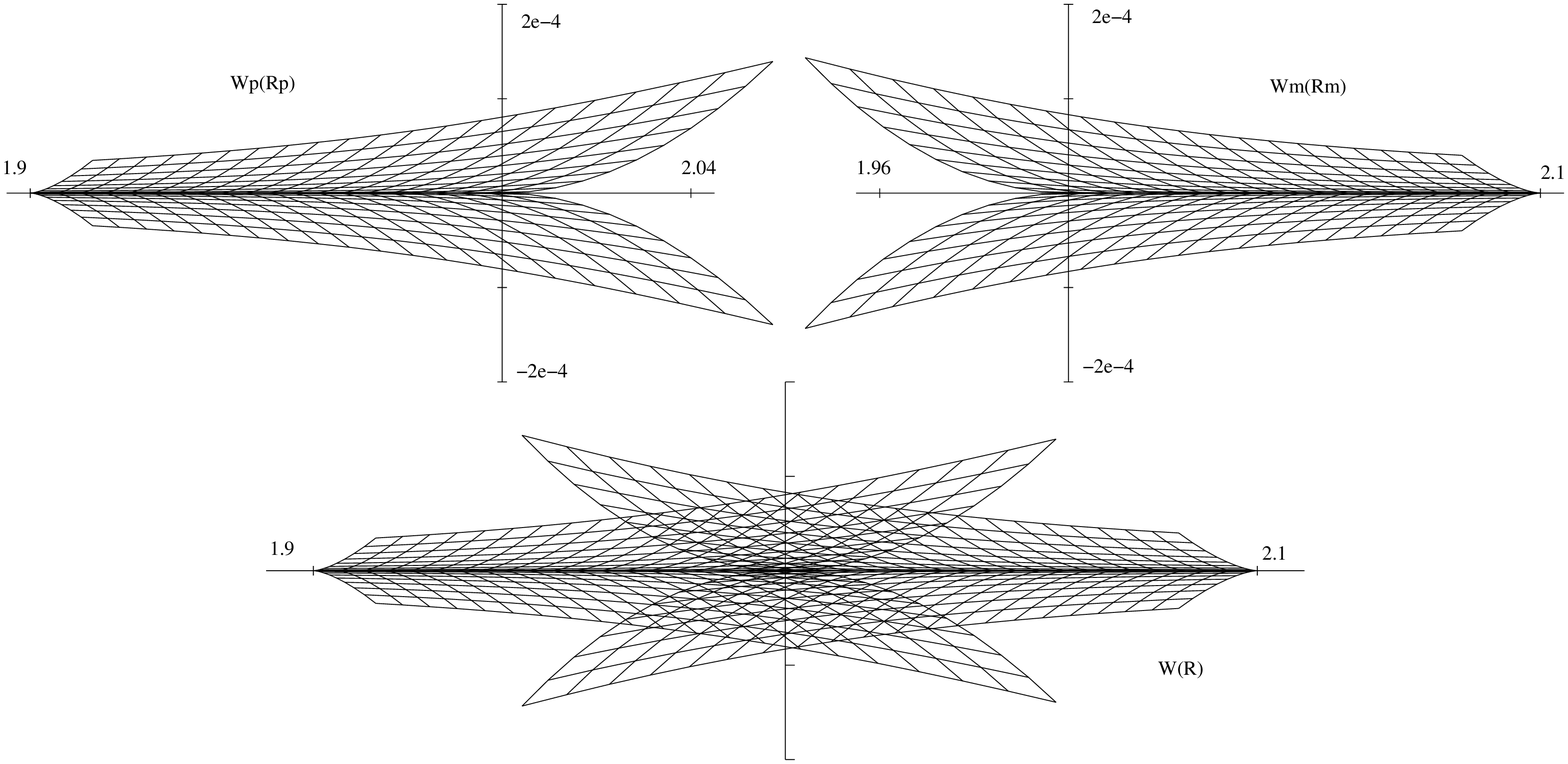}
\end{center}
\caption{$W_+(\Rr_{a,+})$, $W_-(\Rr_{a,-})$ and $W(\Rr_a)$, $a = 1/10$;
vertical scales are stretched.}
\label{fig:wrpm}
\end{figure}

The rectangles $\Rr_a = [2-a,2+a] \times [-a/10,a/10]$ for $a \le 1/10$
are invariant under $W$;
furthermore, the maps $W_\pm$ are injective on $\Rr_{a,\pm}$.
Figure \ref{fig:wrpm} provides strong evidence to these facts,
proved in Propositions \ref{prop:w} and \ref{prop:wdiff}.

\begin{prop}
\label{prop:w}
Take $a \le 1/10$ satisfying also $\Rr_a \subseteq \Uu_3$.
Then $W(\Rr_a) \subseteq \Rr_a$.
Moreover, $|Y(x,y)| \le |y|/49$ for all $(x,y) \in \Rr_a$,
$X_+(x,y) \ge x$ for all $(x,y) \in \Rr_{a,+}$ and
$X_-(x,y) \le x$ for all $(x,y) \in \Rr_{a,-}$.
Also, $W_\pm(\{2\mp a\} \times [-\ver,\ver]) \subset \Rr_{a,\pm}$,
$W_\pm(r_v) \subset \Rr_{a,\mp}$.
\end{prop}

\proof
From Lemma \ref{lemma:omega}, we have $|\omega_\pm(x,y)| \le 2 \ver \le 1/50$
and therefore $|Y(x,y)| \le |y|/49$.
We now prove that $W_+(\Rr_{a,+}) \subset \Rr_a$.
Clearly, $(x,y) \in \Rr_{a,+}$ implies $|Y(x,y)| \le \ver$
so we must prove that $2-a \le X_+(x,y) \le 2+a$.
From equation (\ref{eq:W}),
$X_+(x,y) \ge x$ with equality exactly when $y = 0$.
Since $(\omega_+)_x (x,y) \ge 0$ we have $X_x(x,y) \ge 0$ and
it therefore suffices to prove the two claims in the statement, i.e.,
$X_+(2-a,y) < 2$ and $X_+(2,y) < 2+a$.
For $a \le 1/10$, the inequalities follow from
$\omega_+(x,y) \le 2\ver < a/(4+a) < a/(4-a)$.
Similar checks apply to $W_-$.
\qed

Each rectangle $\Rr_{a,\pm}$ is not invariant.
Denote the interior of a set $X \subset \RR^2$ by $\interior(X)$.
Given $b \ge 0$, $b < a$, let $\Dd_{b,+}$ (resp. $\Dd_{b,-}$)
be the triangle defined by $2-a \le x \le 2-b-10|y|$
(resp. $2+b+10|y| \le x \le 2+a$).
Notice that $\interior(\Xx_0) =
\Rr_a \smallsetminus (\Dd_{0,+} \cup \Dd_{0,-})$.
Let $\Yy_0 \subset \Xx_0$ be the thinner open wedge
defined by $|y| > 10|x-2|$.

\begin{figure}[ht]
\psfrag{x0y0}{$(x_0,y_0)$}
\psfrag{x1y1}{$(x_1,y_1)$}
\psfrag{2-b}{$2-b$}
\psfrag{2-a}{$2-a$}
\psfrag{2+a}{$2+a$}
\psfrag{Db+}{$\Dd_{b,+}$}
\psfrag{X0}{$\Xx_0$}
\psfrag{Y0}{$\Yy_0$}
\psfrag{rh}{$r_h$}
\psfrag{rv}{$r_v$}
\begin{center}
\epsfig{height=35mm,file=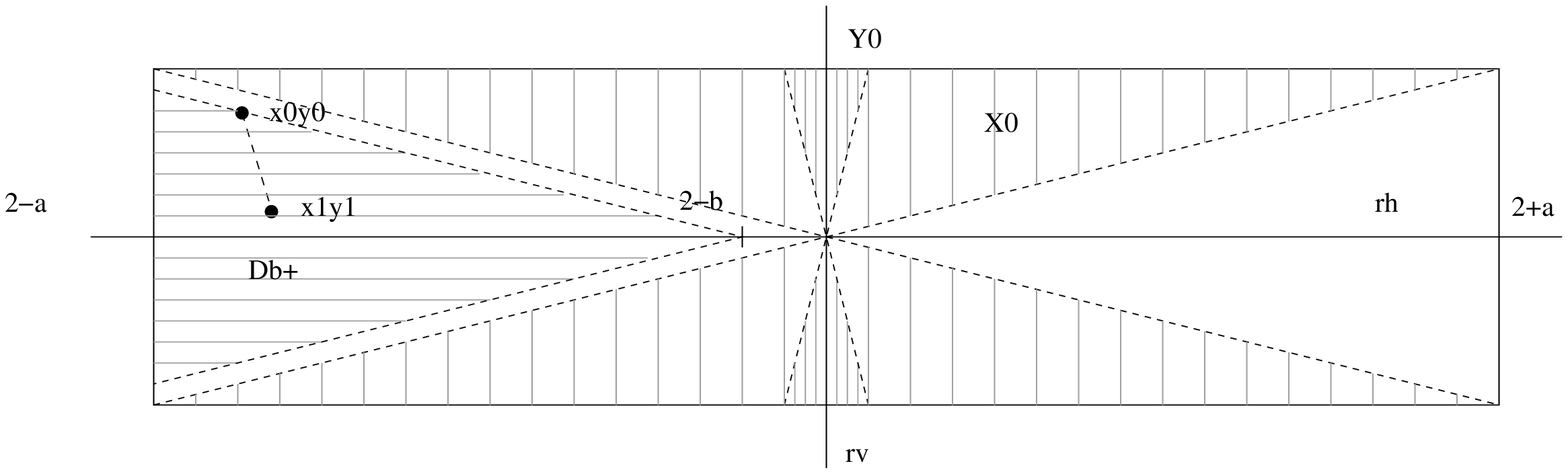}
\end{center}
\caption{$(x_0,y_0) \in \Dd_{b,+}$ implies
$(x_1,y_1) = W_+(x_0,y_0) \in \Dd_{b,+}$.}
\label{fig:tab1}
\end{figure}

\begin{prop}
\label{prop:Tab}
The triangles $\Dd_{b,\pm}$ are invariant, i.e.,
$W_+(\Dd_{b,+}) \subseteq \Dd_{b,+}$ and
$W_-(\Dd_{b,-}) \subseteq \Dd_{b,-}$.
If $p \ne (2\pm b,0), (2\pm a,0)$
is a boundary point of $\Dd_{b,\pm}$
then $W_\pm(p) \in \interior(\Dd_{b,\pm})$.

Finally, $W_+(\Yy_0) \subset \interior(\Dd_{0,-})$
and $W_-(\Yy_0) \subset \interior(\Dd_{0,+})$.
\end{prop}

This will imply (Proposition \ref{prop:rates})
that points in $\Dd_{0,\pm}$ and $\Yy_0$ have cubic convergence
to points in $r_h$ different from $p_\xxxx$.

\proof
We prove that
$W_+(\Dd_{b,+}) \subseteq \Dd_{b,+}$
by computing the slope $\alpha$ of the segment joining $(x_0,y_0)$ and
\[ (x_1,y_1) = W_+(x_0,y_0) = \left(
\frac{1+ \omega_+(x_0,y_0)}{1-\omega_+(x_0,y_0)}\;x_0,
\frac{\omega_+ (x_0,y_0)}{1 + \omega_+(x_0,y_0)}\;y_0 \right), \]
given by
\[ \alpha = \frac{-y_0}{\omega_+} \frac{1 - \omega_+}{2x_0(1+\omega_+)} \]
where $\omega_+$ stands for $\omega_+(x_0,y_0)$
(see Figure \ref{fig:tab1}).
By  Lemma \ref{lemma:omega}, $|\omega_\pm| \le 2|y|$:
simple algebra then shows that, for $a \le 1/10$,
we have $|\alpha| > 1/10$ and the segment is steeper than the
non-vertical sides of $\Dd_{b,+}$.
Since $y_0$ and $y_1$ have the same sign, invariance of $\Dd_{b,+}$ follows.
The argument for $\Dd_{b,-}$ is similar.

The verification that $W_\pm(\Yy_0) \subset \interior(\Dd_{0,\mp})$
uses estimates of the form $|y|/2 < |\omega_\pm(p)| < 1/20$
in the closure of $\Yy_0$; details are left to the reader.
\qed

\begin{prop}
\label{prop:wdiff}
Each map $W_\pm : \interior(\Rr_{a,\pm}) \smallsetminus r_h \to \Rr_a$
is an orientation preserving diffeomorphism to its image.
\end{prop}

\smallskip

\proof
We consider the Jacobian matrix
\begin{equation}
D W_\pm(x,y) = \begin{pmatrix}
\displaystyle\frac{2 (\omega_\pm)_x}{(1 - (\omega_\pm))^2}x +
\frac{1 + (\omega_\pm)}{1 - (\omega_\pm)} &
\displaystyle\frac{2 (\omega_\pm)_y}{(1 - (\omega_\pm))^2}x \\ \\
\displaystyle\pm\frac{(\omega_\pm)_x}{(1 + (\omega_\pm))^2}y &
\displaystyle\pm\frac{(\omega_\pm)_y}{(1 + (\omega_\pm))^2}y \pm
\frac{(\omega_\pm)}{1 + (\omega_\pm)} \end{pmatrix}.
\label{eq:DW}
\end{equation}
From Lemma \ref{lemma:partialomega}, $(X_\pm)_x > 0$ and $(Y_\pm)_y \ge 0$
with equality precisely for $y = 0$.
Similarly, for $y \ne 0$,
$\sign((X_\pm)_y) = \sign((Y_\pm)_x) = \pm \sign(y)$.

To prove that $W_\pm$ are local diffeomorphisms
in $\interior(\Rr_{a,\pm}) \smallsetminus r_h$, write
\[ \det D W_\pm(x,y) =
\pm\frac{1}{1 - \omega_\pm^2}
\left( \frac{2 (\omega_\pm)_x \omega_\pm x}{1 - \omega_\pm} +
(\omega_\pm)_y y + \omega_\pm (1 + \omega_\pm) \right). \]
From Lemma \ref{lemma:partialomega},
all terms in the sum in parenthesis have the same sign
and thus $\det D W_\pm(x,y) > 0$ if $y \ne 0$.

We now prove the injectivity of $W_+$
on $\interior(\Rr_{a,+}) \smallsetminus r_h$.
From symmetry, it suffices to prove the injectivity of $W_+$ on
$\Rr_{a,++} = \{ (x,y) \in \interior(\Rr_{a,+}) \;|\; y \ge 0 \}$.
In other words, given $(x_1,y_1) \in \Rr_a$, $y_1 > 0$,
we must prove that there exists at most one point $(x,y) \in \Rr_{a,++}$
with $W_+(x,y) = (x_1,y_1)$.
Let $\gamma: [0,1] \to \RR^2$ be the piecewise affine counterclockwise
parametrization of the boundary of $\Rr_{a,++}$ with
$\gamma(0) = \gamma(1) = (2-a,0)$, $\gamma(1/4) = (2,0)$,
$\gamma(1/2) = (2,\ver)$ and $\gamma(3/4) = (2-a,\ver)$.
Since $\det DW_+(x,y) > 0$ for $y > 0$ and $y_1 > 0$,
$(x_1,y_1)$ is a regular value of $W_+$
and, assuming $(x_1,y_1) \notin (W_+ \circ \gamma)([0,1])$,
the number of solutions of $W_+(x,y) = (x_1,y_1)$
is given by the winding number $c_1$ of $W_+ \circ \gamma$
around $(x_1,y_1)$.
Recall that a simple way to compute $2c_1$
is to count with signs the intersections of $W_+ \circ \gamma$
with the vertical line through $(x_1,y_1)$.
From the signs of the entries of $DW_+$, the $x$ coordinate
of $(W_+ \circ \gamma)(t)$ is strictly increasing from $t = 0$ to $t = 1/2$
and strictly decreasing from $t = 1/2$ to $t = 1$.
Thus, there are at most two intersection points and $|c_1| \le 1$
implying injectivity.
If $(x_1,y_1) \in (W_+ \circ \gamma)([0,1])$,
the argument above applies to nearby points
and the result follows by continuation outside $r_h$.
The proof of the analogous statement for $W_-$ is similar.
\qed


Actually, the maps $W_\pm: \Rr_{a,\pm} \to \Rr_a$ are homeomorphisms
to their respective images; we omit details.

We conclude the section with some estimates on $DW_\pm$
which will be used in the last section.
A vector $v = (v_1,v_2)$ is \textit{near-horizontal} if
$v_1 > 0$ and $|v_2|/v_1 < 1/25$.

\begin{lemma}
\label{lemma:DW}
Take $p \in \Rr_a$, $p \ne p_\xxxx$.
Let $v = (v_1, v_2)$ be a near-horizontal vector.
Then $\tilde v = DW_\pm(p) v = (\tilde v_1, \tilde v_2)$
is also near-horizontal and $0 < v_1/2 < \tilde v_1 < 10 v_1$.
\end{lemma}


\proof
From the expressions of the entries of
\[ DW_\pm =
\begin{pmatrix} m_{11} & m_{12} \\ m_{21} & m_{22} \end{pmatrix} \]
in equation (\ref{eq:DW}) above
and the estimates in Lemma \ref{lemma:omega}
and Lemma \ref{lemma:partialomega},
\[ 0.96 < m_{11} < 5.5, \quad |m_{12}| < 10.5,
\quad |m_{21}| < 0.0105, \quad |m_{22}| < 0.045. \]
The claims now follow from easy computations.
\qed



\section{Convergence rates of Wilkinson's shift strategy}

We now consider the asymptotic behavior of Wilkinson's shift strategy.
Given a $\bW$-orbit $(T_k)$, $T_{k+1} = \bW(T_k)$,
we study the decay of the entry $(T_k)_{23}$.
In $(x,y)$ coordinates, the $W$-orbits $(p_k)$ are defined by
$p_{k+1} = W(p_k)$, $p_0 \in \Rr_a$;
the relevant issue is the convergence rate of the $y$ coordinate since,
from Proposition \ref{prop:phi},
$(\phi(x,y))_{32}/y = 2r_1/r_2^2 > 0$
is bounded and bounded away from $0$ in $\Rr_a$.

A $W$-orbit $(p_k)$ has
\textit{strictly quadratic (resp. cubic) convergence}
if there exist constants $c, C >0$ such that, for all $k \in \NN$,
\[ c |y_k|^r \le |y_{k+1}| \le C |y_k|^r \]
for $r = 2$ (resp. $r = 3$) (here $p_k = (x_k,y_k)$).



\begin{prop}
\label{prop:rates}
Let $p \in \Rr_a \smallsetminus r_h$, $p_k = W^k(p)$.
If $p_k \in \Dd_{0,\pm}$ for some $k$ then
the $W$-orbit $(p_k)$ has strictly cubic convergence.
Otherwise $p_k \in \Xx_0$ for all $k$ and
convergence is strictly quadratic.
\end{prop}

\proof
The case $y_k = 0$ is trivial.
Assume without loss that $p_k \in \Dd_{0,+}$.
From Proposition \ref{prop:Tab},
there exists $b > 0$ such that
$p_{k+1} \in \Dd_{b,+}$ 
and therefore $p_j \in \Dd_{b,+}$ for all $j > k$.
From Proposition \ref{prop:omega},
$\omega_+$ is an even, smooth function in the $y$-variable
in $\Dd_{b,+}$.
From compactness and Lemma \ref{lemma:omega},
given $b > 0$ there exists $C > 0$ such that
for all $(x,y) \in \Dd_{b,+}$ we have
$y^2 \le \omega_+(x,y) \le C y^2$.
From equation (\ref{eq:W}),
\[ |y|^3/2 \le |Y_+(x,y)| = \frac{\omega_+}{1 + \omega_+} |y| \le 2C |y|^3 \]
for $(x,y) \in \Dd_{b,+}$, yielding strictly cubic convergence.

We now consider orbits in the wedge $\Xx_0 = \Xx_{0,+} \cup \Xx_{0,-}$.
From Lemma \ref{lemma:omega},
$|y|/5 \le |\omega_\pm(x,y)| \le 2|y|$
for $(x,y) \in \Xx_{0,\pm}$:
strictly quadratic convergence now follows from
\[ |y|^2/10 \le |Y_+(x,y)| =
\frac{\omega(x,y)}{2 + \omega(x,y)} |y| \le 4 |y|^2, \]
completing the proof.
\qed

Notice that the constant $C$ in the proof depends on $b$
and therefore the rate of cubic convergence is not uniform in $\Dd_{0,+}$,
consistent with strictly quadratic convergence for orbits in $\Xx_0$.

\begin{coro}
\label{coro:X}
Given $p \in \Rr_a \smallsetminus r_h$,
consider the $W$-orbit $p_k = W^k(p)$.
The following conditions are equivalent:
\begin{enumerate}[(a)]
\item{$\lim_{k \to \infty} p_k = p_\xxxx$;}
\item{$p_k \in \Xx_0$ for all $k$;}
\item{the $W$-orbit $(p_k)$ has strictly quadratic convergence.}
\end{enumerate}
\end{coro}

\proof
The estimate $|Y(x,y)| \le |y|/50$ guarantees convergence
to some point $p_\infty$ in $r_h$.
Orbits contained in $\Xx_0$ must then converge to $p_\xxxx$.
Conversely, if $p_\infty \ne p_\xxxx$ then
$p_\infty \in \interior(\Dd_{b,\pm})$ for some $b > 0$
and $p_k \in \Dd_{b,\pm}$ for sufficiently large $k$.
\qed

Let $\Xx \subset \Rr_a$ be the set of points $p$
for which $\lim p_k = p_\xxxx$.
From Proposition \ref{prop:Tab},
$\Xx \subset \Xx_0 \smallsetminus \Yy_0$ and
therefore $r_v \cap \Xx = r_h \cap \Xx = \{p_\xxxx\}$.
We still need to prove that $\Xx \ne \{p_\xxxx\}$:
this and more will be done in this section.

\begin{figure}[ht]
\begin{center}
\psfrag{p0}{$p_\xxxx$}
\psfrag{30}{$(3,0)$}
\epsfig{height=30mm,file=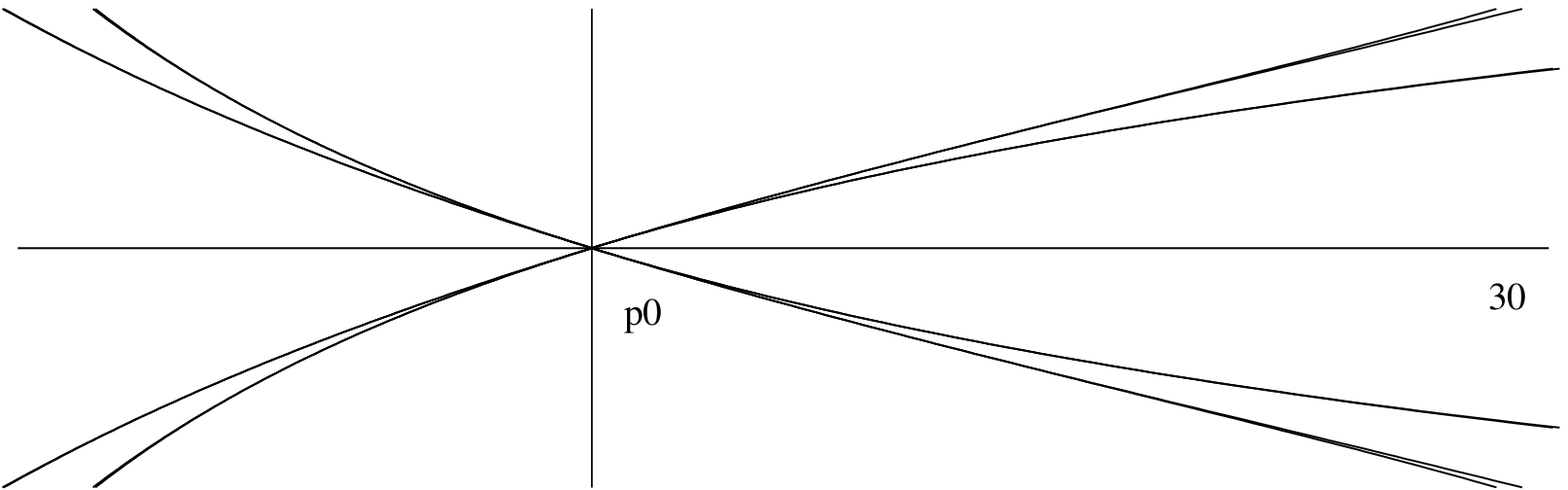}
\end{center}
\caption{\capsize The set $\Xx$ near $p_\xxxx$; in scale.}
\label{fig:xxx}
\end{figure}

Figure \ref{fig:xxx} shows $\Xx$ extended
to a rectangle much larger than $\Rr_{1/10}$:
numerical evidence indicates that even in such larger regions
the qualitative descriptions remain valid.

A compact set $K \subset \RR$ is a \textit{Cantor set}
if $K$ has empty interior and no isolated points.
As we shall prove in Theorem \ref{theo:hausdorff},
horizontal sections of $\Xx$ are Cantor sets.
The set $\Xx$ is the union of \textit{uncountably many} arcs,
disjoint except at $p_\xxxx$.
Each arc intersects a horizontal line in a single point.

The Hausdorff dimension of the middle-third Cantor
is $\log 2/\log 3 \approx 0.63$
(\cite{Falconer} and \cite{Federer} contain
a thorough discussion of Hausdorff dimension).
More generally,
self-similar Cantor sets have positive Hausdorff dimension.
The horizontal sections of $\Xx$ are much thinner:
they have Hausdorff dimension $0$.
That is why the fine structure is invisible in this figure,
unlike most figures of Cantor sets in the quoted books.

Numerical evidence also indicates that the northwest-southeast
leg of set $\Xx$ is the union of a family of analytic curves,
tangent (not crossing) at $p_\xxxx$:
we shall not pursue this matter further.

\begin{figure}[ht]
\begin{center}
\psfrag{X0}{$\Xx_0$}
\psfrag{X1p}{$\Xx_1^{(+)}$}
\psfrag{X1m}{$\Xx_1^{(-)}$}
\psfrag{X2pp}{$\Xx_2^{(++)}$}
\psfrag{X2pm}{$\Xx_2^{(+-)}$}
\psfrag{X2mp}{$\Xx_2^{(-+)}$}
\psfrag{X2mm}{$\Xx_2^{(--)}$}
\epsfig{height=25mm,file=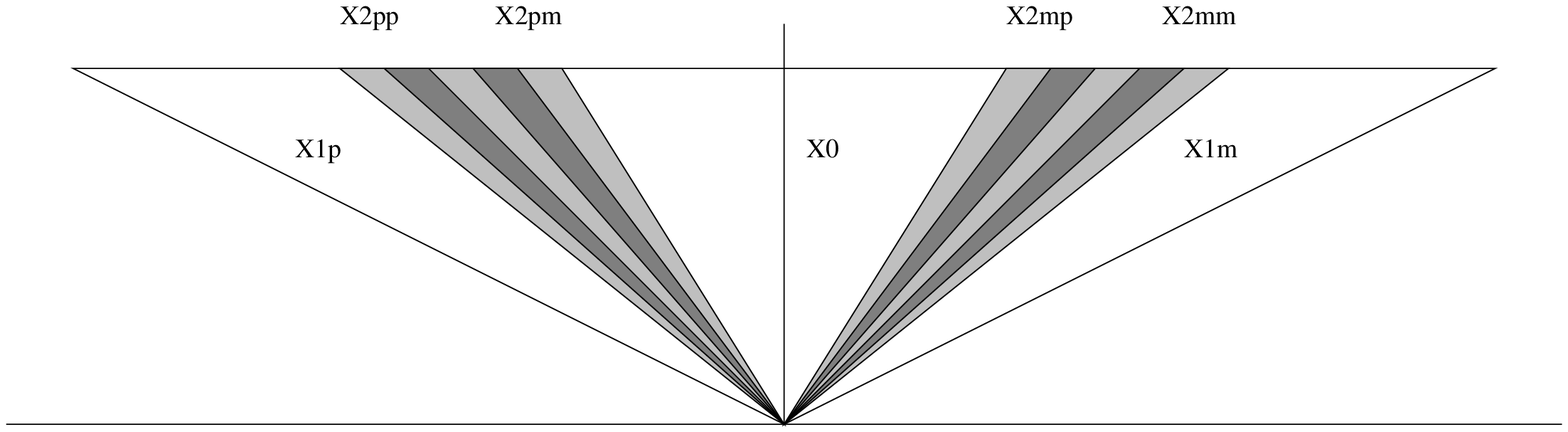}
\end{center}
\caption{\capsize The upper halves of the sets $\Xx_0$,
$\Xx_1^{\pm}$ and $\Xx_2^{\pm\pm}$; schematic.}
\label{fig:x0x1x2}
\end{figure}

We provide another description of $\Xx$ as the intersection
of a nested sequence of compact sets $\Xx_n$.
From Corollary \ref{coro:X},
$p = p_0 \in \Xx$ if and only if $p_k = W^k(p) \in \Xx_0$ for all $k$.
Define $\Xx_n \subset \Rr_a$ to be the set of points $p \in \Rr_a$
such that $p_k \in \Xx_0$ for $0 \le k \le n$.
Thus $\Xx_0 = \Xx_0$, $\Xx_{n+1} = W^{-1}(\Xx_n) \subset \Xx_n$
and $\Xx = \bigcap_n \Xx_n$.
Figure \ref{fig:x0x1x2} indicates the first few sets.
As the diagram suggests, the interior of $\Xx_n$ has $2^{n+1}$
connected components which we now describe.
In a sense, the Cantor sets on horizontal lines
are limits of this successive doubling of components.

\begin{figure}[ht]
\begin{center}
\psfrag{A}{$A$}
\psfrag{B}{$B$}
\psfrag{C}{$C$}
\psfrag{W+-(A)}{$W_+^{-1}(A)$}
\psfrag{W+-(B)}{$W_+^{-1}(B)$}
\psfrag{W+-(C)}{$W_+^{-1}(C)$}
\psfrag{X1}{$\Xx_1^{(+)}$}
\psfrag{pdia}{$p_\xxxx$}
\psfrag{X0capW+}{$\Xx_0 \cap W_+(\Rr_{a,+})$}
\psfrag{W+Ra+}{$W_+(\Rr_{a,+})$}
\epsfig{height=38mm,file=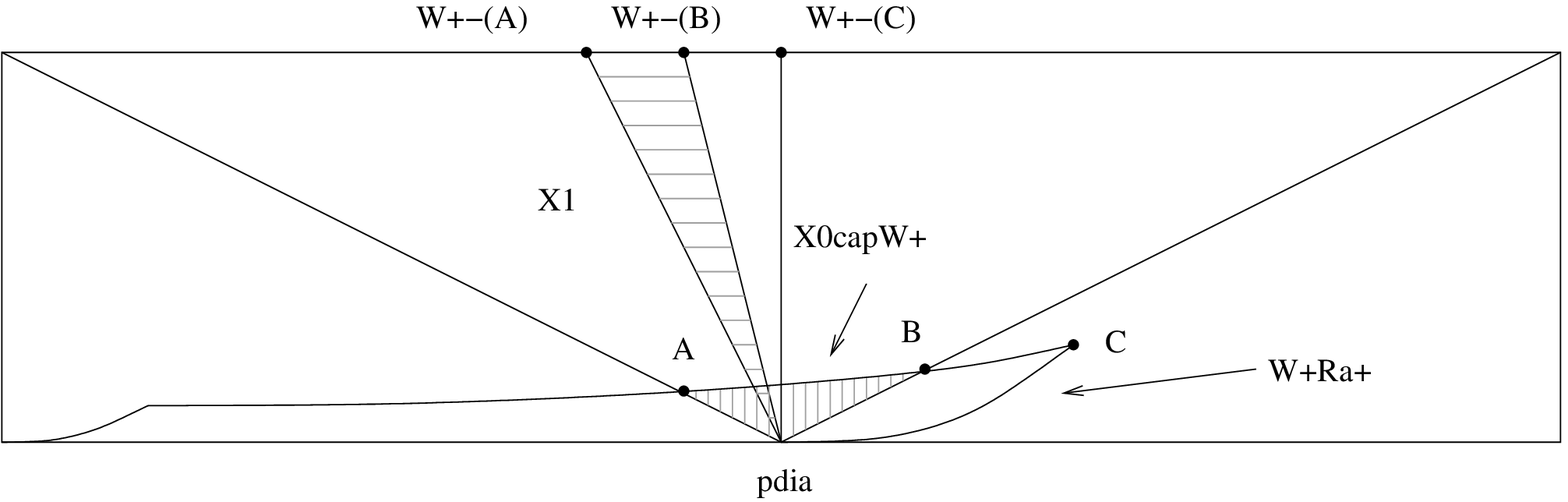}
\end{center}
\caption{\capsize The upper halves of $\Xx_0 \cap W_+(\Rr_{a,+})$
and $\Xx_1^{(+)} = W_+^{-1}(\Xx_0)$; schematic.}
\label{fig:X1}
\end{figure}

Decompose the sets $\Xx_n$
by tracking on which side of $r_v$ the points $p_k$ lie.
A \textit{sign sequence} of length $n$ is
a function $\tau: \{0,\ldots,n-1\} \to \{+,-\}$
or, equivalently, a string of $n$ signs
$(\tau(0), \tau(1), \ldots, \tau(n-1))$.
Define
\[ \Xx_n^{\tau}
= \{ p \in \Rr_a \;|\; p_k \in \Xx_{0,\tau(k)}, 0 \le k < n;
p_n \in \Xx_0 \}. \]
Figure \ref{fig:X1} shows a simple example: the upper half of $\Xx_1^{(+)}$
is a curvilinear triangle with base contained in the top side of $\Rr_{a,+}$
and vertex in $p_\xxxx$;
the upper half of $\Xx_1^{(-)}$ is a similar triangle in $\Rr_{a,-}$.
Also, $\Xx_1^{(+)} = W_+^{-1}(\Xx_0) \cap \Rr_{a,+} \subset \Xx_0$.

We also consider \textit{infinite sign sequences},
i.e., functions $\sigma: \NN \to \{+,-\}$,
and the corresponding sets
\[ \Xx^\sigma = \bigcap_{n \in \NN} \Xx_n^{\sigma|_{\{0,1,\ldots,n-1\}}}. \]
Let $\Sigma$ be the set of infinite sign sequences:
$\Sigma$ admits the natural metric
$d(\sigma_1,\sigma_2) = 2^{-n}$ where $n$ is the smallest
number for which $\sigma_1(n) \ne \sigma_2(n)$.
The metric space $\Sigma$ is homeomorphic to the
middle-third Cantor set contained in $[0,1]$.
The sets $\Xx^\sigma$ are compact, the intersection of two distinct
such sets is $\{p_\xxxx\}$ and the union of all $\Xx^\sigma$ is $\Xx$.
We conclude this section by proving
that the intersection of each $\Xx^\sigma$
with a horizontal line is not empty.

Let $\ell_\pm$ be the sides $y = \pm a^2$, $-a \le x \le a$, of $\Rr_a$.

\begin{prop}
\label{prop:Xsigma}
Let $\gamma_h: [0,1] \to \Rr_a \smallsetminus \{p_\xxxx\}$
be a parametrized curve
with $\gamma_h(0) \in \Dd_{0,+}$, $\gamma_h(1) \in \Dd_{0,-}$.
Let $\sigma$ be an infinite sign sequence:
there exists $t \in [0,1]$ such that $\gamma_h(t) \in \Xx^\sigma$.
\end{prop}

Notice that we do not claim that such $t$ is unique:
this requires stronger hypothesis and will be discussed
in the next section.

\proof
We first prove by induction on $n$
that the connected component containing $p_\xxxx$ of
each set $\Xx_n^\tau$ has elements on the sides $\ell_\pm$
(here $\tau$ is a sign sequence of length $n$).
The case $n = 0$ is trivial
and the case $n = 1$ has already been discussed (Figure \ref{fig:X1}).
Assume the connected component of $\Xx_n^\tau$ containing $p_\xxxx$
to have a point in $\ell_+$.
Consider a path $\gamma_v: [0,1] \to \Xx_n^\tau$
with $\gamma_v(0) = p_\xxxx$, $\gamma_v(1) \in \ell_+$.
The image of $\gamma_v$ must intersect the curve $W_+(\ell_+)$:
let $t_0$ be the smallest $t$ in this intersection.
The image of the restriction of
$\gamma_v$ to $[0,t_0]$ is contained in $W_+(\Rr_{a,+})$.
Since $W_+$ is a homeomorphism on its image (Proposition \ref{prop:wdiff}),
there exists $\gamma_{v+}:[0,t_0] \to \Rr_{a,+}$
with $W_+(\gamma_{v+}(t)) = \gamma_v(t)$ for all $t$.
Also, the image of $\gamma_{v+}$ is contained in $\Xx_{n+1}^{(+,\tau)}$.
A similar construction works for $\ell_-$ and for $\Xx_{n+1}^{(-,\tau)}$,
completing the proof of the claim.

Consider an infinite sign sequence $\sigma$
and its restrictions $\tau_n = \sigma|_{\{0,1,\ldots,n-1\}}$. 
For each $n$ let 
\[ K_n = \{ t \in [0,1] \;|\; \gamma_h(t) \in
\Xx_n^{\tau_n} \}. \]
The sets $K_n$ are nested, compact and nonempty
and therefore their intersection is also nonempty.
Any $t$ in the intersection satisfies $\gamma_h(t) \in \Xx^\sigma$.
\qed

As we shall see in the last section,
the sets $\Xx_n^\tau$ and $\Xx^\sigma$
are connected but the proof of this fact requires more careful estimates.

\section{Geometry of $\Xx$}

The main result of this section, Theorem \ref{theo:hausdorff},
is that $\Xx$ is very thin, almost as thin as a finite union of curves.
More precisely, $\Xx$ has Hausdorff dimension 1.

Denote the length of an interval $I$ by $\mu(I)$.
Given an infinite sign sequence $\sigma$,
define $\sigma^\sharp$ by $\sigma^\sharp(k) = \sigma(k+1)$
(this is the standard shift operator in symbolic dynamics).
Similarly, for a sign sequence $\tau$ of length $n+1$
let $\tau^\sharp$ be the sign sequence of length $n$
defined by $\tau^\sharp(k) = \tau(k+1)$.

\begin{theo}
\label{theo:hausdorff}
For any infinite sign sequence $\sigma$, the set $\Xx^\sigma$
is a curve $\Xx^\sigma = \{ (g^\sigma(y),y), y \in [-a/10,a/10]\}$
where $g^\sigma: [-a/10,a/10] \to [2-a,2+a]$ is a Lipschitz function.
Wilkinson's step takes curves to curves:
$W(\Xx^{\sigma}) \subset \Xx^{\sigma^\sharp}$.
The set $\Xx$ has Hausdorff dimension $1$
and the intersection of $\Xx$ with any horizontal line
is a Cantor set of Hausdorff dimension $0$.
\end{theo}

We need a few preliminary definitions.
A {\it near-horizontal curve} is a $C^1$ function
$\gamma: I \to \Rr_a \smallsetminus \{p_\xxxx\}$
such that, for all $t$ in the interval $I$,
the tangent vector $\gamma'(t)$ is near-horizontal
(as defined at the end of Section 4).
A near-horizontal curve $\tilde\gamma$ is \textit{standard}
if the first coordinate of $\tilde\gamma(t)$ is $x = t$.
For any near-horizontal curve $\gamma$
there exists a unique strictly increasing $C^1$ function $\alpha$,
the \textit{standard reparametrization} of $\gamma$,
for which $\tilde\gamma = \gamma \circ \alpha$ is standard.
The \textit{height} of a near-horizontal curve $\gamma$
is $y_\ast = \tilde\gamma(2)$ so that $\gamma$
crosses $r_v$ at $(2,y_\ast)$.

Near-horizontality is preserved by Wilkinson's step:
from Lemma \ref{lemma:DW},
if $\gamma: I \to \Rr_{a,+} \smallsetminus \{p_\xxxx\}$
(resp. $\Rr_{a,-} \smallsetminus \{p_\xxxx\}$)
is a near-horizontal curve
then so is $W_+ \circ \gamma$ (resp. $W_- \circ \gamma$).
This process squeezes near-horizontal curves towards the line $r_h$.
The constant $1/25$ in the definition is somewhat arbitrary
but it can not be replaced by very small numbers
since $W$ does not take horizontal lines to horizontal lines.

\begin{figure}[ht]
\begin{center}
\psfrag{t0}{$t_0$}
\psfrag{t1}{$t_1$}
\psfrag{gammat(t0)}{$\tilde\gamma(t_0)$}
\psfrag{gammat(t1)}{$\tilde\gamma(t_1)$}
\psfrag{gammat(2)}{$\tilde\gamma(2) = (2,y_\ast)$}
\psfrag{Ta0+}{$\Dd_{0,+}$}
\psfrag{Ta0-}{$\Dd_{0,-}$}
\psfrag{pdia=20}{$p_\xxxx = (2,0)$}
\psfrag{gamma}{$\gamma$}
\epsfig{height=40mm,file=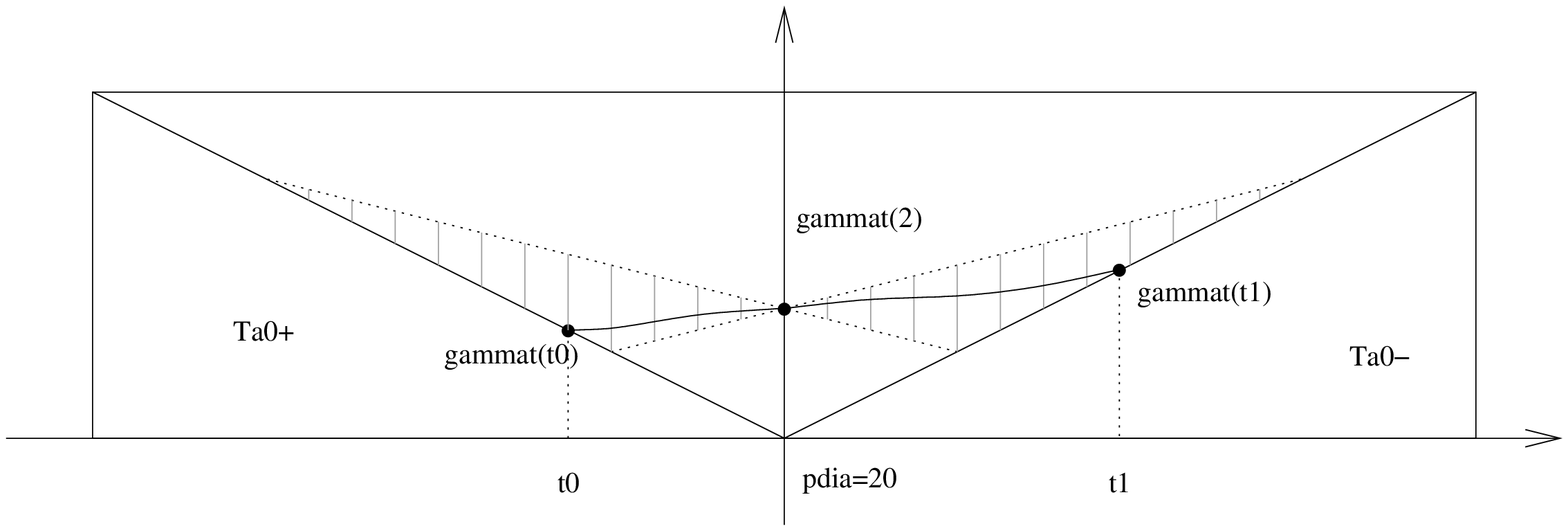}
\end{center}
\caption{\capsize A near-horizontal curve (not in scale).}
\label{fig:step0}
\end{figure}

Recall that $\Yy_0 \subset \Xx_0$ is an open wedge with the property
that $p \in \Yy_0$ implies $W_\pm(p) \in \Dd_{0,\mp}$.
For a sign sequence $\tau$ of length $n$, define
\[ \Yy_n^\tau
= \{ p \in \Rr_a \;|\; p_k \in \Xx_{0,\tau(k)}, 0 \le k < n;
p_n \in \Yy_0 \}. \]
In particular, the orbit $(p_k)$ with $p_0 \in \Yy_n^\tau$
escapes $\Xx_0$ starting from $k = n+1$.
Also, $\Yy_n^\tau$ is an open subset of $\Xx_n^\tau$
disjoint from $\Xx_{n+1}$.

\begin{lemma}
\label{lemma:sizeXtau}
Let $\tilde\gamma: [t_0,t_1] \to \Rr_a \smallsetminus \{p_\xxxx\}$
be a standard near-horizontal curve.
Assume that $\tilde\gamma(t_0) \in \Dd_{0,+}$,
$\tilde\gamma(t_1) \in \Dd_{0,-}$.
Let $y_\ast = \tilde\gamma(2)$ be the height of $\tilde\gamma$.
Then, for $\tilde\gamma(t) = (t,y) \in \Xx_0$,
\[ {|y_\ast|}/{3} < |y| < 3|y_\ast|. \]
Let $\tau$ be a sign sequence of length $n$. The sets
\[ I_n^\tau =
\{ t \in [t_0,t_1] \;|\; \tilde\gamma(t) \in \Xx_n^\tau \}, \quad
J_n^\tau =
\{ t \in [t_0,t_1] \;|\; \tilde\gamma(t) \in \Yy_n^\tau \} \]
are intervals and their lengths satisfy
\[ \frac{1}{10^{n+1}} \left| \frac{y_\ast}{90} \right|^{2^n}
< \mu(J_n^\tau) < \mu(I_n^\tau) <
2^n \left| 40 y_\ast \right|^{2^n} < 2^{(-2^{(n-1)})}. \]
\end{lemma}

\proof
As a basis for an inductive proof,
we first consider the case $n = 0$.
Set
\[ I_0 = \{ t \in [t_0,t_1] \;|\; \tilde\gamma(t) \in \Xx_0 \}, \quad
J_0 = \{ t \in [t_0,t_1] \;|\; \tilde\gamma(t) \in \Yy_0 \}. \]
Draw lines through the point $(2,y_\ast)$
with linear coefficients $\pm 1/25$,
as in Figure \ref{fig:step0},
and compute their intersections with the diagonals of $\Rr_a$.
Since these diagonals are steeper than $\tilde\gamma$,
the intersection of the image of $\tilde\gamma$ with $\Xx_0$
is contained in the shaded triangles.
An elementary geometric argument proves the first claim,
verifies that $I_0$ and $J_0$ are intervals
and obtains the estimates for their lengths.

\begin{figure}[ht]
\begin{center}
\psfrag{Xn}{$\Xx_n^{\tau^\sharp}$}
\psfrag{Xn+1}{$\Xx_{n+1}^{\tau}$}
\psfrag{yn}{$y_\ast^\sharp$}
\psfrag{yn+1}{$y_\ast$}
\psfrag{gamman}{$\tilde\gamma^\sharp$}
\psfrag{gamman+1}{$\tilde\gamma$}
\psfrag{In}{$I_n^{\tau^\sharp}$}
\psfrag{In+1}{$I_{n+1}^{\tau}$}
\psfrag{Ta0+}{$\Dd_{0,+}$}
\psfrag{Ta0-}{$\Dd_{0,-}$}
\psfrag{xhyh}{$(\hat x, \hat y)$}
\epsfig{height=40mm,file=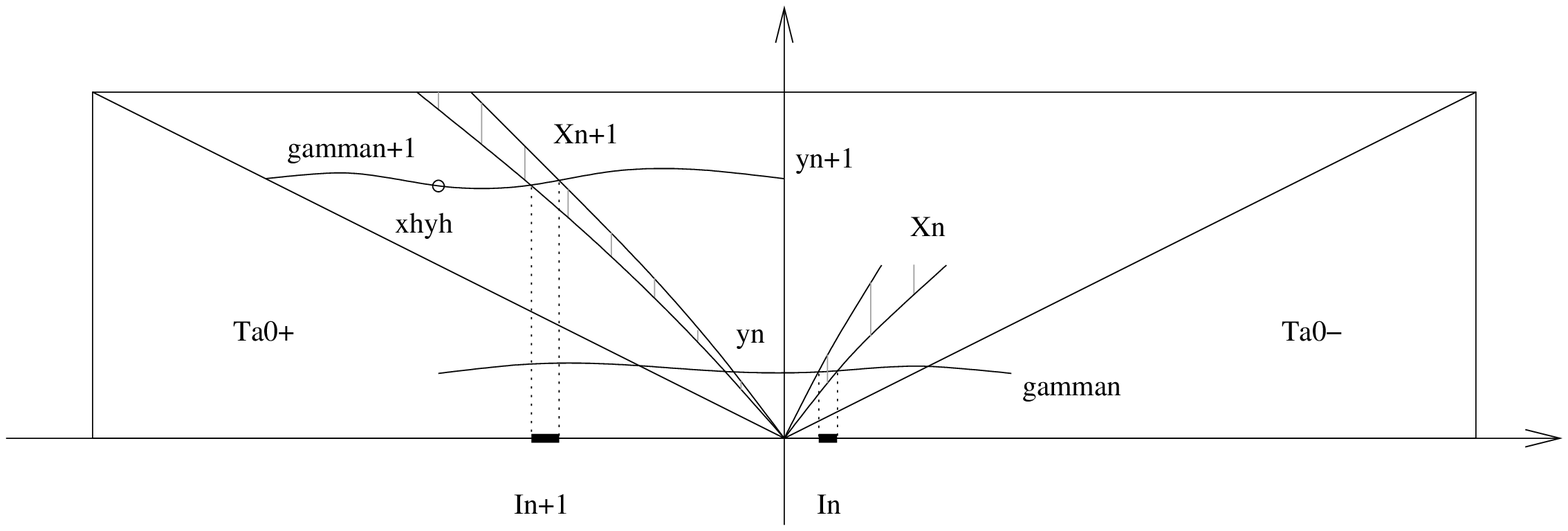}
\end{center}
\caption{\capsize A near-horizontal curve (not in scale).}
\label{fig:stepn}
\end{figure}

We now do the induction step.
For a sign sequence $\tau$ of length $n+1$
with $\tau(0) = +$ (the other case is similar),
consider
\[ I_{n+1}^\tau =
\{ t \in [t_0,t_1] \;|\; \tilde\gamma(t) \in \Xx_{n+1}^\tau \} \]
and $\gamma_+$, the restriction of $\tilde\gamma$ to $[t_0,2]$.
Set $\gamma^\sharp = W_+ \circ \gamma_+$ (see Figure \ref{fig:stepn}):
as remarked above, $\gamma^\sharp$ is a near-horizontal curve.
Let $\alpha^\sharp: [t_0^\sharp, t_1^\sharp] \to [t_0,2]$
be the standard reparametrization of $\gamma^\sharp$ so that
$\tilde\gamma^\sharp = \gamma^\sharp \circ \alpha^\sharp$
is a standard near-horizontal curve with height
$y_\ast^\sharp = \tilde\gamma^\sharp(2)$.
Notice that $\tilde\gamma^\sharp(t_0^\sharp) \in \Dd_{0,+}$
and $\tilde\gamma^\sharp(t_1^\sharp) \in \Dd_{0,-}$
so, by the induction hypothesis, the sets
\[ I_n^{\tau^\sharp} =
\{ t \in [t_0^\sharp,t_1^\sharp] \;|\;
\tilde\gamma^\sharp(t) \in \Xx_n^{\tau^\sharp} \}, \quad
J_n^{\tau^\sharp} =
\{ t \in [t_0^\sharp,t_1^\sharp] \;|\;
\tilde\gamma^\sharp(t) \in \Yy_n^{\tau^\sharp} \} \]
are intervals whose lengths $\ell^\sharp_I, \ell^\sharp_J$ satisfy
\begin{equation}
\frac{1}{10^{n+1}} \left| \frac{y_\ast^\sharp}{90} \right|^{2^n}
< \ell^\sharp_J < \ell^\sharp_I <
2^n \left| 40 y_\ast^\sharp \right|^{2^n}.
\label{eq:stepn}
\end{equation}
Clearly, $t \in I_n^{\tau^\sharp}$ (resp. $J_n^{\tau^\sharp}$)
if and only if
$\alpha^\sharp(t) \in I_{n+1}^{\tau}$ (resp. $J_{n+1}^{\tau}$),
proving that $I_{n+1}^{\tau}$ and $J_{n+1}^{\tau}$ are intervals;
let $\ell_I, \ell_J$ be their lengths.
By Lemma \ref{lemma:DW}, $1/10 < (\alpha^\sharp)'(t) < 2$
for all $t \in I_n$ and therefore
\begin{equation}
(1/10) \ell^\sharp_J < \ell_J < \ell_I < 2 \ell^\sharp_I.
\label{eq:elln}
\end{equation}
Let $\tilde\gamma(\alpha^\sharp(2)) = (\hat x, \hat y)$
so that $W_+(\hat x,\hat y) = (2,y_\ast^\sharp)$.
From the first claim, $|y_\ast|/3 < |\hat y| < 3|y_\ast|$.
By Lemma \ref{lemma:omega},
$(1/5) |\hat y| < \omega_+(\hat x,\hat y) < 2|\hat y|$.
Thus, by equation \ref{eq:W},
$(1/10) \hat y^2 < |y_\ast^\sharp| < 3 \hat y^2$.
Combining these estimates,
\begin{equation}
\left| \frac{y_\ast}{90} \right| > \left|\frac{y_\ast^\sharp}{90} \right|^2, \quad
\left| 40 y_\ast^\sharp \right| < \left| 40 y_\ast \right|^2.
\label{eq:squares}
\end{equation}
The required estimate now follows from estimates
\ref{eq:stepn}, \ref{eq:elln}, \ref{eq:squares}
and the fact that $|y_\ast| \le a/10 \le 1/100$.
\qed

{\nobf Proof of Theorem \ref{theo:hausdorff}: }
We organize the proof by stating a few claims,
most of which in terms of a standard near-horizontal curve $\tilde\gamma$;
notation will be borrowed from Lemma \ref{lemma:sizeXtau}.
For a sign sequence $\tau$ of length $n$,
consider the sign sequences $(\tau,+)$ and $(\tau,-)$
of length $n+1$.

\smallbreak
\noindent\textbf{(a)}
\textit{The intervals $I_{n+1}^{(\tau,+)}$ and $I_{n+1}^{(\tau,-)}$
are contained in different connected components
of $I_n^{\tau} \smallsetminus J_n^{\tau}$.}

This claim is proved by induction.
The case $n = 0$ follows from the inclusions
$\Xx_1^\pm \subset \Xx_{0,\pm} \smallsetminus \Yy_0$.
The induction step employs the same construction used in
Lemma \ref{lemma:sizeXtau}: for
$\tilde\gamma^\sharp = W^+ \circ \tilde\gamma \circ \alpha^\sharp$,
the standard reparametrization $\alpha^\sharp$ maps the intervals
constructed from $\tau^\sharp$ to the corresponding intervals
for $\tau$.

\smallbreak
\noindent\textbf{(b)}
\textit{For an infinite sign sequence $\sigma$
there is a unique $t^\sigma \in [t_0,t_1]$
with $\tilde\gamma(t^\sigma) \in \Xx^\sigma$.}

This follows directly from Lemma \ref{lemma:sizeXtau}
since 
\[ \lim_{n \to \infty} \mu(I_n^{\sigma|_{\{0, 1, \ldots, n-1\}}}) = 0. \]
In particular, a horizontal line at height $y$
(which is near-horizontal) meets $\Xx^\sigma$
at a single point $(g^\sigma(y),y)$.
The Lipschitz estimate
$|g^\sigma(y_1) - g^\sigma(y_2)| \le 25 |y_1 - y_2|$
likewise follows by taking lines of slope less than $1/25$.

\smallbreak
\noindent\textbf{(c)}
\textit{Let $\sigma_1$ and $\sigma_2$ be distinct sign sequences
and $n$ be the smallest integer for which $\sigma_1(n) \ne \sigma_2(n)$.
Then
\[ \frac{1}{10^{n+1}} \left| \frac{y_\ast}{90} \right|^{2^n}
< |t^{\sigma_1} - t^{\sigma_2}| <
2^n \left| 40 y_\ast \right|^{2^n} < 2^{(-2^{(n-1)})}. \]
}

Let $\tau$ be the restriction
of $\sigma_1$ or $\sigma_2$ to $\{0,1,\ldots, n-1\}$.
Without loss,
$t^{\sigma_1} \in I_{n+1}^{(\tau,+)}$ and
$t^{\sigma_2} \in I_{n+1}^{(\tau,-)}$:
Lemma \ref{lemma:sizeXtau} and claim (a) obtain the inequalities.

In particular, the map from the Cantor set $\Sigma$ to
$\Kk = \{t \in [t_0,t_1]\;|\;\tilde\gamma(t) \in \Xx\}$
taking $\sigma$ to $t_\sigma$ is a homeomorphism.

We now construct a very thin open cover of $\Kk$.
For a positive integer $n$
let $S^+_n \subset \Sigma$ be the set of the $2^n$ sign sequences
$\sigma$ for which $k \ge n$ implies $\sigma(k) = +$.
Let $B_r(p)$ be the ball of radius $r$ and center $p$.

\smallbreak
\noindent\textbf{(d)}
\textit{Take $r > 0$ and let $n$ be a positive integer
such that $r > 2^{(- 2^{(n-1)})}$.
Then the balls $B_r(t^{\sigma^+}) \subset \RR$, $\sigma^+ \in S^+_n$,
form an open cover of $\Kk$.}

From step (b), given $t \in \Kk$
there is some $\sigma \in \Sigma$ for which $t = t^\sigma$.
Let $\sigma_n \in S^+_n$ be defined by
\[ \sigma_n(k) = \begin{cases}
\sigma(k), &0 \le k < n, \\ +, &k \ge n \end{cases} \]
so that $\sigma_n \in S^+_n$.
From claim (c), $|t^{\sigma_n} - t^{\sigma}| < 2^{(- 2^{(n-1)})}$
and therefore $t \in B_r(t^{\sigma_n})$.

The open cover of $\Kk$ above is now adapted to yield
a thin open cover of $\Xx$.
For a positive integer $m$, let $Y_m \subset [-a/10,a/10]$
be a set of $2m$ equally spaced points so that given $y \in [-a/10,a/10]$
there is $y_0 \in Y_m$ with $|y-y_m| < a/(10 m)$.
Also, for positive integers $n$ and $m$, let $X_{n,m} \subset \Xx$
be the set of $m 2^{(n+1)}$ elements of the form
$(g^\sigma(y),y)$, $\sigma \in S^+_n$, $y \in Y_m$.

\smallbreak
\noindent\textbf{(e)}
\textit{Consider $r > 0$ and positive integers $n$ and $m$ such that
$r > 2^{(1 - 2^{(n-1)})}$ and $r > 1/m$.
Then the balls $B_r(p) \subset \RR^2$, $p \in X_{n,m}$,
form an open cover of $\Xx$.}

Given $p = (g^\sigma(y),y)$ let $y_0$ be the point of $Y_m$
closest to $y$ and let $\sigma_n \in S^+_n$ be defined as in claim (d).
Let $p_1 = (g^{\sigma_n}(y), y)$ and $p_0 = (g^{\sigma_n}(y_0), y_0)$.
Notice that $p_0 \in X_{n,m}$.
In order to prove that $p$ belongs to $B_r(p_0)$
we estimate $|p-p_1|$ and $|p_1-p_0|$.
From claim (d), $|p-p_1| < r/2$.
Since $a \le 1/10$, $|y - y_0| < a/(10m) \le r/100$;
from the Lipschitz estimate after claim (b),
$|g^{\sigma_0}(y) - g^{\sigma_0}(y_0)| \le r/4$
and therefore $|p_1 - p_0| < r/2$.

\smallbreak
\noindent\textbf{(f)}
\textit{The set $\Kk$ has Hausdorff dimension $0$.}

Given $r > 0$, let $n = \lceil 1 + \log_2(-\log_2 r) \rceil$.
Claim (d) obtains on open cover of $\Kk$ by
$2^n < -4\log_2 r$ balls of radius $r$.
For all $d > 0$ we have
\[ \lim_{r \to 0} (-4\log_2 r) r^d = 0. \]
Thus, $\mu_d(\Kk)$,
the Hausdorff measure of dimension $d$ of $\Kk$, equals $0$.

\smallbreak
\noindent\textbf{(g)}
\textit{The set $\Xx$ has Hausdorff dimension $1$.}

Since $\Xx$ contains curves this dimension is at least one.
It suffices therefore to prove that $\mu_d(\Xx) = 0$ for all $d > 1$.

Given $r > 0$, let $n = \lceil 2 + \log_2(-\log_2 r) \rceil$
and $m = \lceil1/r\rceil$.
Claim (e) obtains an open cover of $\Xx$ by
$m 2^{n+1} < -32(\log_2 r) r^{-1}$ balls of radius $r$.
For all $d > 1$ we have
\[ \lim_{r \to 0} -32(\log_2 r) r^{-1} r^d = 0. \]
This implies $\mu_d(\Xx) = 0$, proving the claim
and completing the proof of the theorem.
\qed

\section{Conclusion}

Recall that $\bomega_+(T)$ and $\bomega-(T)$ are the roots of 
a quadratic polynomial with smooth coefficients (in $T$).
The algorithm to obtain $\bomega$ given $T$ is thus multivalued:
in variables $x$ and $y$, the domains of choice of each
sign are rectangles sharing the line $r_v$.
In a sense, dynamical systems techniques arose naturally in this paper
because orbits may switch sides arbitrarily often and
according to arbitrary patterns (sign sequences).
An outstanding example of the use of sign sequences to label orbits
are kneading sequences (\cite{MilnorThurston}).
Many numerical algorithms do not exhibit this kind of complexity.

Unsurprisingly, the set $\Xx$ of atypical points
is more complicated than an algebraic surface.
This is why Hausdorff dimension becomes the relevant tool
to quantify the size of $\Xx$.
Theorem \ref{theo:hausdorff} shows that $\Xx$ is very thin.

A  matrix is an {\it AP-matrix} if its spectrum contains an
arithmetic progression with three terms and is {\it AP-free} otherwise.
This paper concentrated on the simplest examples of AP-matrices,
with eigenvalues $-1$, $0$ and $1$.
We conjecture that the presence of three term arithmetic progressions
in the spectrum (and whether the three eigenvalues are consecutive)
dictates the possible rates of convergence of Wilkinson's shift.
In a forthcoming paper we show that for AP-free tridiagonal matrices,
cubic convergence of Wilkinson's shift indeed holds
(\cite{LST2}, \cite{LST3}).

%

\smallskip

{\noindent\bf Acknowledgements: }
The editorial board of FoCM did an extraodinary job in articulating
a diverse collection of referees which led to a far more readable text.
We are very grateful to these editors and referees.
The authors acknowledge support from CNPq, CAPES, IM-AGIMB and Faperj.

\smallbreak
{

\parindent=0pt
\parskip=0pt
\obeylines

Ricardo S. Leite, Departamento de Matemática, UFES
Av. Fernando Ferrari, 514, Vitória, ES 29075-910, Brazil

\smallskip

Nicolau C. Saldanha and Carlos Tomei, Departamento de Matem\'atica, PUC-Rio
R. Marqu\^es de S. Vicente 225, Rio de Janeiro, RJ 22453-900, Brazil

\smallskip

rsleite@pq.cnpq.br
nicolau@mat.puc-rio.br; http://www.mat.puc-rio.br/$\sim$nicolau/
tomei@mat.puc-rio.br

}


\begin{thebibliography}{[10]}

\bibitem{BS}{ Batterson, S. and Smillie, J.,
{Rayleigh quotient iteration for nonsymmetric matrices},
Math. of Comp., 55, 169-178, 1990.}
\bibitem{BG}{ de Boor, C. and Golub, G. H.,
{The numerically stable reconstruction of a Jacobi matrix from spectral data},
Lin. Alg. Appl. 21, 245-260, 1978.}
\bibitem{DNT}{ Deift, P., Nanda, T., Tomei, C.,
{Differential equations for the symmetric eigenvalue problem},
SIAM J. Num. Anal. 20, 1-22, 1983.}
\bibitem{Demmel}{ Demmel, J. W.,
{Applied Numerical Linear Algebra},
SIAM, Philadelphia, 1997.}
\bibitem{Falconer}{ Falconer, K. J.,
{Fractal geometry: mathematical foundations and applications},
John Willey, Chichester, 2003.}
\bibitem{Federer}{ Federer, H.,
{Geometric measure theory},
Classics in Mathematics, Springer-Verlag, New York, 1996.}
\bibitem{GH}{ Gragg, W. B. and Harrod, W. J.,
{The numerically stable reconstruction of Jacobi matrices from spectral data},
Numer. Math. 44, 317-335, 1984.}
\bibitem{LST1}{ Leite, R. S., Saldanha, N.C. and Tomei, C.,
{An atlas for tridiagonal isospectral manifolds},
Linear Algebra Appl. 429, 387-402, 2008.}
\bibitem{LST2}{ Leite, R. S., Saldanha, N.C. and Tomei, C.,
{The asymptotic of Wilkinson's shift:
cubic convergence for generic spectra}, in preparation.}
\bibitem{LST3}{ Leite, R. S., Saldanha, N.C. and Tomei, C.,
{The asymptotics of Wilkinson's shift iteration},
preprint, {\tt http://www.arxiv.org/abs/math.NA/0412493}.}
\bibitem{MilnorThurston}{ Milnor, J. and Thurston, W.,
{On iterated maps of the interval}, in
Dynamical Systems (ed. Alexander, J. C.),
Lecture Notes in Math. 1342, 465-563, New York: Springer-Verlag, 1988.}
\bibitem{Moser}{ Moser, J.,
{Finitely many mass points on the line under the influence
of an exponential potential},
in Dynamic systems theory and applications (ed. J.~Moser),
Lecture Notes in Phys. 38, 467-497, New York, 1975.}
\bibitem{Parlett}{ Parlett, B. N.,
{The Symmetric Eigenvalue Problem},
Prentice-Hall, Englewood Cliffs, NJ 1980.}
\bibitem{Tomei}{ Tomei, C.,
{The Topology of Manifolds of Isospectral Tridiagonal Matrices},
Duke Math. J., 51, 981-996, 1984.}
\bibitem{Wilkinson}{ Wilkinson, J. H.,
{The algebraic eigenvalue problem},
Oxford University Press, 1965.}


\end{thebibliography}
\end{document}